\documentclass[11pt,a4paper]{article}
\usepackage{amsmath}
\usepackage{multirow}
\usepackage{slashbox}

\def\inter{\mathop{\hbox{\rm int}}}

\def\HypR{H_{{\cal R}}}
\def\HypS{H_{\cal S}}
\def\HypRs{H_{{\cal R}i}}
\def\HypSs{H_{{\cal S}i}}

\usepackage[bf,normalsize,tableposition=top]{caption}
\usepackage{subfig}
\usepackage{amsfonts}
\usepackage{color}
\usepackage{epsfig}
\usepackage{srcltx}
\usepackage[colorlinks=true]{hyperref}
\usepackage{amssymb}
\def\D{{\cal D}}
\def\Z{{\cal Z}}

\def\Erf{{\mathop{\hbox{\small\rm Erf}}}}

\newcommand{\half}{ \mbox{\small$\frac{1}{2}$}}

\newcommand{\eight}{ \mbox{\small$\frac{1}{8}$}}
\def\SadVal{\mathop{\hbox{\rm SadVal}}}

\usepackage{graphicx}
\usepackage{algorithm}
\usepackage[colorlinks=true]{hyperref}
\def\bR{{\mathbf{R}}}
\def\bZ{{\mathbf{Z}}}
\def\Opt{{\hbox{\rm Opt}}}

\def\A{{\cal A}}

\def\U{{\cal U}}
\def\SG{{\cal SG}}
\def\T{{\cal T}}
\def\E{{\cal E}}
\def\P{{\cal P}}
\def\Q{{\cal Q}}
\def\R{{\cal R}}
\newcommand{\I}{{\cal I}}
\def\G{{\cal G}}
\def\H{{\cal H}}
\def\M{{\cal M}}
\def\X{{\cal X}}

\def\N{{\cal N}}
\def\C{{\cal C}}
\def\S{{\cal S}}

\def\bE{{\mathbf{E}}}
\def\bS{{\mathbf{S}}}
\def\Prob{\hbox{\rm Prob}}
\def\Erf{\hbox{\rm Erf}}
\newtheorem{lemma}{Lemma}[section]
\newtheorem{proposition}{Proposition}[section]

\def\Det{{\hbox{\rm Det}}}
\def\qed{$\Box$}

\def\e{{\rm e}}
\def\Tr{{\hbox{\rm Tr}}}

\definecolor{MyDarkBlue}{rgb}{0,0.08,0.45}
\definecolor{MyViolet}{rgb}{0.45,0.08,0.95}
\definecolor{MyBrown}{rgb}{0.45,0.08,0}

\newcommand{\be}{\begin{eqnarray}}
\newcommand{\ee}[1]{\label{#1}\end{eqnarray}}
\newcommand{\nn}{\nonumber \\}
\newcommand{\ese}{\end{eqnarray*}}
\newcommand{\bse}{\begin{eqnarray*}}
\newcommand{\rf}[1]{~(\ref{#1})}

\title{Hypothesis Testing via Affine Detectors}
\author{
Anatoli Juditsky
\thanks{LJK, Universit\'e Grenoble Alpes, B.P. 53, 38041 Grenoble Cedex 9, France,	
{\tt anatoli.juditsky@imag.fr}}
\and Arkadi Nemirovski
\thanks{Georgia Institute
 of Technology, Atlanta, Georgia
30332, USA, {\tt nemirovs@isye.gatech.edu}\newline
The first author was supported by the CNRS-Mastodons project GARGANTUA,
and the LabEx PERSYVAL-Lab (ANR-11-LABX-0025). Research of
the second author was supported by NSF grants  CMMI-1262063, CCF-1523768.}}
\date{}
\begin{document}
\maketitle
\begin{abstract}
In this paper, we further develop the approach, originating in \cite{GJN}, to ``computation-friendly'' hypothesis testing
via Convex Programming. Most of the existing results on hypothesis testing aim to quantify in a closed analytic form separation between sets of distributions allowing for reliable decision in precisely stated observation models. In contrast to this descriptive (and highly instructive) traditional framework, the approach
we promote here can be qualified as operational -- the testing routines and their risks are yielded by an efficient computation. All we know in advance is that, under favorable circumstances, specified in \cite{GJN},  the risk of such test, whether high or low, is provably near-optimal under the circumstances.
As a compensation for the lack of ``explanatory power,'' this approach
is applicable to a much wider family of observation schemes and hypotheses to be tested than those where ``closed form descriptive analysis'' is possible.
\par
In the present paper our primary emphasis is on computation: we make a step further in extending the principal tool developed in \cite{GJN} -- testing routines based on affine detectors --  to a large variety of testing problems. The price of this development is the loss of blanket near-optimality of the proposed procedures (though it is still preserved in the observation schemes studied in \cite{GJN}, which now become particular cases of the general setting considered here).
\end{abstract}
\section{Introduction}
 This paper can be considered as an extension of \cite{GJN} where the following simple observation was the starting point of numerous developments:
\begin{quote}Imagine that we want to decide on two composite hypotheses about the distribution $P$ of a random observation $\omega$ taking values in observation space $\Omega$, $i$-th hypothesis stating that $P\in\P_i$, where $\P_i$, $i=1,2$, are given families of probability distributions on $\Omega$. Let $\phi:\;\Omega\to \bR$ be a {\sl detector}, and let the {\sl risk} of detector $\phi$ be defined as the smallest $\epsilon_\star$ such that
\begin{equation}\label{riskriskrisk}
\int_\Omega\e^{-\phi(\omega)}P(d\omega)\leq\epsilon_\star\,\,\forall P\in\P_1\ \&\
\int_\Omega\e^{\phi(\omega)}P(d\omega)\leq\epsilon_\star\,\,\forall P\in\P_2.
\end{equation}
Then the test $\T^K$ which, given $K$ i.i.d. observations $\omega_t\sim P\in\P_1\cup\P_2$, $t=1,...,K$, deduces that $P\in\P_1$ when $\sum_{t=1}^K\phi(\omega_t)\geq0$, and that $P\in\P_2$ otherwise, accepts the true hypothesis with $P$-probability at least $1-\epsilon_\star^K$.
\end{quote}
It was shown in \cite{GJN} that
\begin{enumerate}
\item\label{item1} the detector-based tests are ``near optimal'' -- if the above hypotheses can be decided upon by a single-observation test $\T^*$ with risk $\delta<1/2$, there exists a detector $\phi$ with ``comparable'' risk $\epsilon_\star=2\sqrt{\delta(1-\delta)}$.

    Note that while the risk $2\sqrt{\delta(1-\delta)}$ seems to be much larger than $\delta$, especially for small $\delta$, we can ``compensate'' for risk deterioration passing from the single-observation test $\T^1$ associated with the detector $\phi$ to the test $\T^K$ based on the same detector and using $K$ observations. The risk of the test $\T^K$, by the above, is upper-bounded by $\epsilon_\star^K=[2\sqrt{\delta(1-\delta)}]^K$ and thus is not worse than the risk $\delta$ of the ``ideal'' single-observation test already for a ``quite moderate'' value $\rfloor {2\over 1-\ln(4(1-\delta))/\ln(1/\delta)}\lfloor$ of $K$.
\item There are ``good,'' in certain precise sense, parametric families of distributions, primarily,
\begin{itemize}
\item Gaussian $\N(\mu,I_d)$ distributions on $\Omega=\bR^d$,
\item Poisson distributions with parameters $\mu\in\bR^d_+$ on $\Omega=\bZ^d$; the corresponding random variables have $d$ independent entries, $j$-th of them being Poisson with parameter $\mu_j$,
\item Discrete distributions on $\{1,...,d\}$, the parameter $\mu$ of a distribution being the vector of probabilities to take value $j=1,...,d$,
\end{itemize}
{for which} the optimal (with the minimal risk, and thus -- near-optimal by \ref{item1}) detectors can be found efficiently, provided that $\P_i$, $i=1,2$, are {\sl convex hypotheses}, meaning that they are cut off the
family of distributions in question by {restricting} the distribution's parameter $\mu$ {to reside in a convex domain.}
\footnote{In retrospect, these results can be seen as a development of the line of research initiated by the pioneering works of H. Chernoff \cite{chernoff1952}, C. Kraft \cite{Kraft1955}, and L. Le Cam \cite{Lecam1973}, further developed in \cite{Birge1980,birge1981,Birge1982,Birge1983,Burnashev1979,Burnashev1982} among many others (see also references in \cite{Ingster2002}).}
\end{enumerate}
On a closer inspection, the ``common denominator'' of  Gaussian, Poisson and Discrete families of distributions is that in all these cases the minimal risk detector for a pair of convex hypotheses is {\sl affine},\footnote{affinity of a detector makes sense only when $\Omega$ can be naturally identified with a subset of some $\bR^d$. This indeed is the case for Gaussian and Poisson distributions; to make it the case for discrete distributions on $\{1,...,d\}$, it suffices to encode $j\leq d$ by $j$-th basic orth in $\bR^d$, thus making $\Omega$ the set of basic orths in $\bR^d$. With this encoding, {\sl every} real valued function on $\{1,...,d\}$ becomes affine.} and the results of \cite{GJN} in the case of deciding on a pair of convex hypotheses  stemming from a {\em good family of distributions} sum up to the following:
\begin{itemize}
\item[A)]  the best -- with the smallest possible risk -- {\sl affine} detector, same as its risk, can be efficiently computed;
\item[B)] the smallest risk {\sl affine} detector from A) is the best, in terms of risk, detector available under the circumstances, so that the associated test is near-optimal.
\end{itemize}
Note that as far as practical applications of the above approach are concerned, one ``can survive'' without B) (near-optimality of the constructed detectors), while A) {\em is a must.}
In this paper, we focus on families of distributions obeying A); this class turns out to be incomparably larger than
what was defined as ``good'' in \cite{GJN}. In particular, it includes nonparametric families of distributions. Staying within this much broader class, we still are able to construct in a computationally efficient way the best affine detectors for a pair of ``convex'', in certain precise sense, hypotheses, along with valid upper bounds on the risks of the detectors. What we, in general, can{\sl not} claim anymore, is that the tests associated with the above detectors are near-optimal. This being said, we believe that investigating possibilities for building tests and quantifying their performance in a computationally friendly manner is of value even when we cannot provably guarantee near-optimality of these tests.
\par
The  paper is organized as follows. The families of distributions well suited for constructing affine detectors in a computationally friendly fashion are introduced and investigated in section \ref{sect:Setup}. In particular, we develop a kind of fully algorithmic ``calculus'' of these families. This calculus demonstrates that the families of probability distributions covered by our approach are much more common commodity than ``good observation schemes'' as defined in \cite{GJN}.  In section \ref{sect:testing} we explain how to build within our framework tests for pairs (and larger tuples) of hypotheses and how to quantify performance of these tests in a computationally efficient fashion. Aside of general results of this type, we work out in detail the case where the family of distributions giving rise to ``convex hypotheses'' to be tested is comprised of sub-Gaussian distributions (section \ref{subGaussian}). In section \ref{sect:aggr} we discuss an application to the  now-classical statistical problem -- aggregation of estimators -- and show how the results of \cite{Golden2009} can be extended to the general situation considered here. Finally, in section \ref{sec:quadratic} we show how our framework can be extended in the Gaussian case to include quadratic detectors.
To streamline the presentation, all {proofs exceeding few lines} are collected in the appendix.

\section{Setup}\label{sect:Setup}
Let us fix {\sl observation space} $\Omega=\bR^d$, and let $\P_j$, $1\leq j\leq J$, be given families of Borel probability distributions on $\Omega$. Put broadly, our goal is, given a random observation $\omega\sim P$, where $P\in\bigcup\limits_{j\leq J}\P_j$, to decide upon the hypotheses $H_j:P\in\P_j$, $j=1,...,J$. We intend to address the above goal in the case when the families $\P_j$ are {\sl simple} -- they are comprised of distributions for which moment-generating functions
admit an explicit upper bound.
\subsection{Regular and simple probability distributions}
Let
\begin{itemize}
\item $\H$ be a nonempty closed convex set in $\Omega=\bR^d$ symmetric w.r.t. the origin,
\item $\M$ be a closed convex set in some $\bR^n$,
\item $\Phi(h;\mu):\H\times\M \to\bR$ be a continuous function convex in $h\in \H$ and concave in $\mu\in \M$.
\end{itemize}
We refer to $\H,\M,\Phi(\cdot,\cdot)$ satisfying the above restrictions as to {\sl regular data}. Regular data
$\H,\M,\Phi(\cdot,\cdot)$ define a family
$$
\R=\R[\H,\M,\Phi]
$$
of Borel probability distributions $P$ on $\Omega$ such that
\begin{equation}\label{eq1}
\begin{array}{l}
\forall h\in\H\,\,\exists \mu\in\M: \ln\left(\int_\Omega\exp\{h^T\omega\}P(d\omega)\right)\leq \Phi(h;\mu).
\end{array}
\end{equation}
We say that distributions satisfying \rf{eq1} are {\em regular}, and, given regular data $\H,\M,\Phi(\cdot,\cdot)$, we refer to $\R[\H,\M,\Phi]$ as to {\sl regular} family of distributions associated with the data $\H$, $\M$, $\Phi$.
The same regular data $\H,\M,\Phi(\cdot,\cdot)$ define a smaller family
$$
\S=\S[\H,\M,\Phi]
$$
of Borel probability distributions $P$ on $\Omega$ such that
\begin{equation}\label{eq1a}
\begin{array}{l}
\exists \mu\in\M: \forall h\in\H: \ln\left(\int_\Omega\exp\{h^T\omega\}P(d\omega)\right)\leq \Phi(h;\mu).
\end{array}
\end{equation}
We say that distributions satisfying \rf{eq1a} are {\em simple}. Given regular data $\H,\M,\Phi(\cdot,\cdot)$, we refer to $\S[\H,\M,\Phi]$ as to {\sl simple} family of distributions associated with the data $\H$, $\M$, $\Phi$.
\par
Recall that the starting point of our study is a ``plausibly good'' detector-based test which, given  two families $\P_1$ and $\P_2$ of distribution with common observation space, and independent observations $\omega_1,...,\omega_t$ drawn from a distribution $P\in\P_1\cup \P_2$, decides whether $P\in\P_1$ or $P\in \P_2$.
Our interest in regular/simple families of distributions stems from the fact that when the families $\P_1$ and $\P_2$ are of this type, building such test reduces to solving a convex-concave game and thus can be carried on in a computationally efficient manner. We postpone the related construction and analysis to section \ref{sect:testing}, and continue with presenting some basic examples of simple and regular families of distributions and a simple ``calculus'' of these families.
\subsection{Basic examples of simple families of probability distributions}\label{sect:basicexamples}
\subsubsection{Sub-Gaussian distributions}\label{Example1}  Let $\H=\Omega=\bR^d$, $\M$ be a closed convex subset of the set $\G_d=\{\mu=(\theta,\Theta):\theta\in\bR^d,\Theta\in\bS^d_+\}$, where $\bS^d_+$ is cone of positive semidefinite matrices in the space $\bS^d$ of symmetric $d\times d$ matrices, and let
$$
\Phi(h;\theta,\Theta)=\theta^Th+\half h^T\Theta h.
$$
In this case, $\S[\H,\M,\Phi]$ contains all sub-Gaussian distributions $P$ on $\bR^d$ with sub-Gaussianity parameters from $\M$, that is, all Borel probability distributions $P$ on $\Omega$ admitting upper bound
\begin{equation}\label{eq40}
\bE_{\omega\sim P} \{\exp\{h^T\omega\}\}\leq \exp\{\theta^Th+\half h^T\Theta h\}\;\;\forall h\in\bR^d
\end{equation}
on the moment-generating function, with parameters $(\theta,\Theta)$ of the bound belonging to $\M$. In particular, $\S[\H,\M,\Phi]$ contains all Gaussian distributions $\N(\theta,\Theta)$ with $(\theta,\Theta)\in\M$.
\subsubsection{Poisson distributions}\label{Example2} Let $\H=\Omega=\bR^d$, let $\M$ be a closed convex subset of $d$-dimensional nonnegative orthant  $\bR^d_+$, and let
$$
\Phi(h=[h_1;...;h_d];\mu=[\mu_1;...;\mu_d])=\sum_{i=1}^d\mu_i[\exp\{h_i\}-1]: \H\times \M\to\bR.
$$
\def\cP{\mathop{\hbox{\rm Poisson}}}
The family $\S[\H,\M,\Phi]$ contains all product-type Poisson distributions $\cP[\mu]$ with vectors $\mu$ of parameters belonging to $\M$; here $\cP[\mu]$ is the distribution of random $d$-dimensional vector with independent of each other entries, $i$-th entry being Poisson random variable with parameter $\mu_i$.
\subsubsection{Discrete distributions}\label{Example3} Consider a discrete random variable taking values in $d$-element set $\{1,2,...,d\}$, and let us think of such a variable as of random variable taking values $e_i$, $i=1,...,d$, where $e_i$ are standard basic orths in $\bR^d$; probability distribution of such a variable can be identified with a point $\mu=[\mu_1;...;\mu_d]$ from the $d$-dimensional probabilistic simplex
\def\bDelta{{\mathbf{\Delta}}}
$$
\bDelta_d=\{\nu\in\bR^d_+:\sum_{i=1}^d\nu_i=1\},
$$
where $\mu_i$ is the probability for the variable to take value $e_i$. With these identifications, setting $\H=\bR^d$, specifying $\M$ as a closed convex subset of $\bDelta_d$ and setting
$$
\Phi(h=[h_1;...;h_d];\mu=[\mu_1;...;\mu_d])=\ln\left(\sum_{i=1}^d\mu_i\exp\{h_i\}\right),
$$
the family $\S[\H,\M,\Phi]$ contains distributions of all discrete random variables taking values in $\{1,...,d\}$ with probabilities $\mu_1,...,\mu_d$ comprising a vector from $\M$.
\subsubsection{Distributions with bounded support}\label{sect:boundedsupport} Consider the family $\P[X]$ of Borel probability distributions
supported on a closed and bounded convex set $X\subset \Omega=\bR^d$, and let
$$
\phi_X(h)=\max_{x\in X} h^Tx
$$
be the support function of $X$. We have the following result (to be refined in section \ref{sect:Support}):
\begin{proposition}\label{propbounded} For every $P\in\P[X]$ it holds
\begin{equation}\label{eq090807}
\forall h\in\bR^d: \;\;\ln\left(\int_{\bR^d}\exp\{h^T\omega\}P(d\omega)\right) \leq h^Te[P]+{1\over 8}\left[\phi_X(h)+\phi_X(-h)\right]^2,
\end{equation}
where $e[P]=\int_{\bR^d}\omega P(d\omega)$ {is the expectation of $P$},
and the right hand side function in {\rm (\ref{eq090807})} is convex. As a result, setting
$$
\H=\bR^d,\,\,\M=X,\,\,\Phi(h;\mu)=h^T\mu+{1\over 8}\left[\phi_X(h)+\phi_X(-h)\right]^2,
$$
we get regular data such that $\P[X]\subset \S[\H,\M,\Phi]$.
\end{proposition}
For proof, see Appendix \ref{app1}

\subsection{Calculus of regular and simple families of probability distributions}
Regular and simple families of probability distributions admit ``fully algorithmic'' calculus, with the main calculus rules as follows.
\subsubsection{Direct summation} For $1\leq\ell\leq L$, let regular data $\H_\ell\subset\Omega_\ell=\bR^{d_\ell}$,  $\M_\ell\subset\bR^{n_\ell}$, $\Phi_\ell(h_\ell;\mu_\ell):\H_\ell\times \M_\ell\to\bR$ be given. Let us set
$$
\begin{array}{rcl}
\Omega&=&\Omega_1\times...\times\Omega_L=\bR^{d},\;d=d_1+...+d_L,\\
\H&=&\H_1\times...\times\H_L=\{h=[h^1;...;h^L]:h^\ell\in\H_\ell,\ell\leq L\},\\
\M&=&\M_1\times...\times\M_L=\{\mu=[\mu^1;...;\mu^L]:\mu^\ell\in\M^\ell,\ell\leq L\}\subset\bR^{n},\;n=n_1+...+n_L,\\
\multicolumn{3}{l}{\Phi(h=[h^1;...;h^L];\mu=[\mu^1;...;\mu^L])=\sum_{\ell=1}^L\Phi_\ell(h^\ell;\mu^\ell):\H\times\M\to\bR.}\\
\end{array}
$$
Then $\H$ is a symmetric w.r.t. the origin closed convex set in $\Omega=\bR^d$, $\M$ is  a nonempty closed convex set in $\bR^n$, $\Phi:\;\H\times\M\to\bR$ is a continuous
convex-concave function, and clearly
\begin{itemize}
\item the family $\R[\H,\M,\Phi]$ contains all product-type distributions $P=P_1\times...\times P_L$ on $\Omega=\Omega_1\times...\times\Omega_L$ with $P_\ell\in\R[\H_\ell,\M_\ell,\Phi_\ell]$, $1\leq\ell\leq L$;
\item the family $\S[\H,\M,\Phi]$ contains all product-type distributions $P=P_1\times...\times P_L$ on $\Omega=\Omega_1\times...\times\Omega_L$ with $P_\ell\in\S[\H_\ell,\M_\ell,\Phi_\ell]$, $1\leq\ell\leq L$.
\end{itemize}
\subsubsection{IID summation}
Let $\Omega=\bR^d$ be an observation space, $(\H,\M,\Phi)$ be regular data on this space, and let $\lambda=\{\lambda_\ell\}_{\ell=1}^K$ be a collection of reals. We can associate with the outlined entities a new data
$(\H_\lambda,\M,\Phi_\lambda)$ on $\Omega$ by setting
$$
\H_\lambda=\{h\in\Omega:\,\|\lambda\|_\infty h\in \H\},\;\;\Phi_{\lambda}(h;\mu)=\sum_{\ell=1}^L\Phi(\lambda_\ell h;\mu):\,\H_\lambda\times\M\to\bR.
$$
 Now, given a probability distribution $P$ on $\Omega$, we can associate with it and with the above $\lambda$ a new probability distribution $P^\lambda$ on $\Omega$ as follows: $P^\lambda$ is the distribution of $\sum_\ell\lambda_\ell \omega_\ell$, where $\omega_1,\omega_2,...,\omega_L$  are drawn, independently of each other, from $P$. An immediate observation is that the data $(\H_\lambda,\M,\Phi_\lambda)$ is regular, and
\begin{itemize}
 \item whenever a probability distribution $P$ belongs to $\S[\H,\M,\Phi]$, the distribution $P^\lambda$ belongs to $\S[\H_\lambda,\M,\Phi_\lambda]$. In particular, when $\omega\sim P\in\S[\H,\M,\Phi]$, then the distribution $P^L$ of the sum of $L$ independent copies of $\omega$ belongs to $\S[\H,\M,L\Phi]$.
\end{itemize}
\subsubsection{Semi-direct summation}  For $1\leq\ell\leq L$, let regular data $\H_\ell\subset \Omega_\ell=\bR^{d_\ell}$, $\M_\ell$, $\Phi_\ell$ be given. To avoid complications, we assume that for every $\ell$,
\begin{itemize}
\item $\H_\ell=\Omega_\ell$,
 \item $\M_\ell$ is bounded.
 \end{itemize}
Let also an $\epsilon>0$ be given. We assume that $\epsilon$ is small, namely, $L\epsilon<1$.
\par
Let us aggregate the given regular data into a new one by setting
$$
\H=\Omega:=\Omega_1\times...\times\Omega_L=\bR^{d},\;d=d_1+...+d_L,\;\;\M=\M_1\times...\times\M_L,
$$
and let us define function $\Phi(h;\mu):\Omega^d\times \M\to\bR$ as follows:
\begin{equation}\label{eq30}
\begin{array}{l}
\Phi(h=[h^1;...;h^L];\mu=[\mu^1;...;\mu^L])=\inf_{\lambda\in \bDelta^\epsilon}\sum_{\ell=1}^d\lambda_\ell\Phi_\ell(h^\ell/\lambda_\ell;\mu^\ell),\\
\bDelta^\epsilon=\{\lambda\in\bR^d:\lambda_\ell\geq\epsilon\,\forall \ell\ \&\ \sum_{\ell=1}^L\lambda_\ell=1\}.\\
\end{array}
\end{equation}
By evident reasons, the infimum in the description of $\Phi$ is achieved, and $\Phi$ is continuous. In addition, $\Phi$ is convex in $h\in\bR^d$ and concave in
$\mu\in\M$. Postponing for a moment verification, the consequences are that $\H=\Omega=\bR^d$, $\M$ and $\Phi$ form a regular data. We claim that
\begin{quote}
{\sl Whenever $\omega=[\omega^1;...;\omega^L]$ is a Borel random variable taking values in $\Omega=\bR^{d_1}\times...\times\bR^{d_L}$, and the marginal distributions $P_\ell$, $1\leq \ell\leq L$, of $\omega$ belong to the families $\S[\bR^{d_\ell},\M_\ell,\Phi_\ell]$ for all $1\leq \ell\leq L$, the distribution $P$ of $\omega$ belongs to $\S[\bR^d,\M,\Phi]$.}
\end{quote}
Indeed, since $P_\ell\in \S[\bR^{d_\ell},\M_\ell,\Phi_\ell]$, there exists $\widehat{\mu}^\ell\in\M_\ell$ such that
$$\ln(\bE_{\omega^\ell\sim P_\ell} \{\exp\{g^T\omega^\ell\}\})\leq \Phi_\ell(g;\widehat{\mu}^\ell)\,\,\forall g\in\bR^{d_\ell}.
$$
Let us set $\widehat{\mu}=[\widehat{\mu}^1;...;\widehat{\mu}^L]$, and let $h=[h^1;...;h^L]\in\Omega$ be given. We can find $\lambda\in\bDelta^\epsilon$ such that
$$
\Phi(h;\widehat{\mu})=\sum_{\ell=1}^L\lambda_\ell \Phi_\ell(h^\ell/\lambda_\ell;\widehat{\mu}^\ell).
$$
Applying H\"older inequality, we get
$$
\bE_{[\omega^1;...;\omega^L]\sim P}\left\{\exp\{\sum_\ell [h^\ell]^T\omega^\ell\}\right\}
\leq \prod\limits_{\ell=1}^L \left(\bE_{\omega^\ell\sim P_\ell}\left\{[h^\ell]^T\omega^\ell/\lambda_\ell\right\}\right)^{\lambda_\ell},
$$
whence
$$
\ln\left(\bE_{[\omega^1;...;\omega^L]\sim P}\left\{\exp\{\sum_\ell [h^\ell]^T\omega^\ell\}\right\}\right)\leq\sum_{\ell=1}^L\lambda_\ell\Phi_\ell(h^\ell/\lambda_\ell;\widehat{\mu}^\ell)=\Phi(h;\widehat{\mu}).
$$
We see that
$$
\ln\left(\bE_{[\omega^1;...;\omega^L]\sim P}\left\{\exp\{\sum_\ell [h^\ell]^T\omega^\ell\}\right\}\right)\leq\Phi(h;\widehat{\mu})\,\,\forall h\in\H=\bR^d,
$$
and thus $P\in\S[\bR^d,\M_\ell,\Phi_\ell]$, as claimed.\par
It remains to verify that the function $\Phi$ defined by (\ref{eq30}) indeed is convex in $h\in\bR^d$ and concave in $\mu\in\M$. Concavity in $\mu$ is evident.
Further, functions $\lambda_\ell\Phi_\ell(h^\ell/\lambda_\ell;\mu)$ (as perspective transformation of convex functions $\Phi_\ell(\cdot;\mu)$) are jointly convex in $\lambda$ and $h^\ell$, and so is $\Psi(\lambda,h;\mu)=\sum_{\ell=1}^L \lambda_\ell\Phi_\ell(h^\ell/\lambda_\ell,\mu)$. Thus $\Phi(\cdot ;\mu)$, obtained by partial minimization of $\Psi$ in $\lambda$, indeed is convex.
\subsubsection{Affine image} Let $\H$, $\M$, $\Phi$ be regular data, $\Omega$ be the embedding space of $\H$, and $x\mapsto Ax+a$ be an affine mapping from $\Omega$ to
$\bar{\Omega}=\bR^{\bar{d}}$, and let us set
 $$\bar{\H}=\{\bar{h}\in\bR^{\bar{d}}:A^T\bar{h}\in\H\},\; \bar{\M}=\M,\, \bar{\Phi}(\bar{h};\mu)=\Phi(A^T\bar{h};\mu)+a^T\bar{h}:\,\bar{\H}\times\bar{M}\to\bR.$$
 Note that $\bar{\H}$, $\bar{\M}$ and $\bar{\Phi}$ are regular data. It is immediately seen that
  \begin{quote}
  {\sl Whenever the probability distribution of a random variable $\omega$ belongs to $\R[\H,\M,\Phi]$ (or belongs to $\S[\H,\M,\Phi]$), the distribution $\bar{P}[P]$ of the random variable $\bar{\omega}=A\omega+a$ belongs to  $\R[\bar{\H},\bar{\M},\bar{\Phi}]$ (respectively, belongs to $\S[\bar{\H},\bar{\M},\bar{\Phi}]$)}.
  \end{quote}
 \subsubsection{Incorporating support information}\label{sect:Support}
Consider the situation as follows. We are given regular data $\H\subset\Omega=\bR^d,\M,\Phi$ and are interested in the family of distribution $\P$ known to belong to $\S[\H,\M,\Phi]$. In addition, we know that all distributions $P$ from $\P$ are supported on a given closed convex set $X\subset\bR^d$. How could we incorporate this domain information to pass from the family
$\S[\H,\M,\Phi]$ containing $\P$  to a smaller family of the same type still containing $\P$ ? We are about to give an answer in the simplest case of $\H=\Omega$. Specifically, denoting by $\phi_X(\cdot)$ the support function of $X$
 and selecting somehow a closed convex set $G\subset\bR^d$ containing the origin, let us set
$$
\widehat{\Phi}(h;\mu)=\inf_{g\in G}\left[\Phi^+(h,g;\mu):=\Phi(h-g;\mu)+\phi_X(g)\right],
$$
where $\Phi(h;\mu):\bR^d\times\M\to\bR$ is the continuous convex-concave function participating in the original regular data.
Assuming that $\widehat{\Phi}$ is real-valued and continuous on the domain $\bR^d\times \M$ (which definitely is the case when $G$ is a compact set such that $\phi_X$ is finite and continuous on $G$), note that $\widehat{\Phi}$ is convex-concave on this domain, so that $\bR^d,\M,\widehat{\Phi}$ is a regular data. We claim that
\begin{quote}
{\sl The family $\S[\bR^d,\M,\widehat{\Phi}]$ contains $\P$, provided the family $\S[\bR^d,\M,\Phi]$ does so, and the first of these two families is smaller than the second one.}
\end{quote}
Verification of the claim is immediate. Let $P\in\P$, so that for properly selected $\mu=\mu_P\in\M$ and for all $e\in \bR^d$ it holds
$$
\ln\left(\int_{\bR^d}\exp\{e^T\omega\}P(d\omega)\right)\leq \Phi(e;\mu_P).
$$
Besides this, for every $g\in G$ we have $\phi_X(\omega)-g^T\omega\geq0$ on the support of $P$, whence for every $h\in\bR^d$ one has
$$
\ln\left(\int_{\bR^d}\exp\{h^T\omega\}P(d\omega)\right)\leq \ln\left(\int_{\bR^d}\exp\{h^T\omega+\phi_X(g)-g^T\omega\}P(d\omega)\right)\leq
\phi_X(g)+\Phi(h-g;\mu_P).
$$
Since the resulting inequality holds true for all $g\in G$, we get
$$
\ln\left(\int_{\bR^d}\exp\{h^T\omega\}P(d\omega)\right)\leq\widehat{\Phi}(h;\mu_P)\,\,\forall  h\in \bR^d,
$$
implying that $P\in\S[\bR^d,\M,\widehat{\Phi}]$; since $P\in\P$ is arbitrary, the first part of the claim is justified. The inclusion $\S[\bR^d,\M,\widehat{\Phi}]\subset \S[\bR^d,\M,\Phi]$ is readily given by the inequality $\widehat{\Phi}\leq\Phi$, and the latter is  due to $\widehat{\Phi}(h,\mu)\leq \Phi(h-0,\mu)+\phi_X(0)$.
\paragraph{Illustration: distributions with bounded support revisited.} In section \ref{sect:boundedsupport}, given a convex compact set $X\subset\bR^d$ with support function $\phi_X$, we checked that
the  data $\H=\bR^d$, $\M=X$, $\Phi(h;\mu)=h^T\mu+{1\over 8}[\phi_X(h)+\phi_X(-h)]^2$ are regular and
the family $\S[\bR^d,\M,\Phi]$  contains the family $\P[X]$ of all Borel probability distributions supported on $X$.
Moreover, for every $\mu\in \M=X$, the family $\S[\bR^d,\{\mu\},\Phi\big|_{\bR^d\times\{\mu\}}]$ contains all supported on $X$ distributions with the expectations $e[P]=\mu$. Note that $\Phi(h;e[P])$ describes well the behaviour of the logarithm $F_P(h)=\ln\left(\int_{\bR^d}\e^{h^T\omega}P(d\omega)\right)$ of the  moment-generating function of  $P\in\P[X]$ when $h$ is small  (indeed, $F_P(h)=h^Te[P]+O(\|h\|^2)$ as $h\to 0$),  and by far overestimates $F_P(h)$ when $h$ is large. Utilizing the above construction, we replace $\Phi$ with the real-valued, convex-concave and continuous on $\bR^d\times\M$ function
$$
\widehat{\Phi}(h;\mu)=\inf_g\left[(h-g)^T\mu+{1\over 8}[\phi_X(h-g)+\phi_X(-h+g)]^2+\phi_X(g)\right]\leq\Phi(h;\mu).
$$
It is easy to see that $\widehat{\Phi}(\cdot;\cdot)$  still ensures the inclusion $P\in\S[\bR^d,\{e[P]\},\widehat{\Phi}\big|_{\bR^d\times\{e[P]\}}]$ for every distribution $P\in\P[X]$ and ``reproduces $F_P(h)$ reasonably well'' for both small and large $h$. Indeed, since $F_P(h)\leq \widehat{\Phi}(h;e[P])\leq\Phi(h;e[P])$,
 for small $h$ $\widehat{\Phi}(h;e[P])$ reproduces ${F}_P(h)$ even better than $\Phi(h;e[P])$, and we clearly have
$$
\widehat{\Phi}(h;\mu)\leq \left[(h-h)^T\mu+{1\over 8}[\phi_X(h-h)+\phi_X(-h+h)]^2+\phi_X(h)\right]=\phi_X(h)\,\,\forall \mu,
$$
and $\phi_X(h)$ is a correct description  of $F_P(h)$ for large $h$.
\section{Affine detectors and hypothesis testing}\label{sect:testing}
\subsection{Situation}\label{sect:situ}
Assume we are given two collections of regular data with common $\Omega=\bR^d$ and $\H$, specifically, the collections  $(\H,\M_\chi,\Phi_\chi)$, $\chi=1,2$.
We start with a construction of a specific test for a pair of hypotheses $H_1:P\in\P_1$, $H_2:P\in\P_2$, where
$$
\P_\chi=\R[\H,\M_\chi,\Phi_\chi],\,\,\chi=1,2.
$$
When building the test, we impose on the regular data in question the following
\begin{quote}
{\bf Assumption I:} {\sl The regular data  $(\H,\M_\chi,\Phi_\chi)$, $\chi=1,2$, are such that the convex-concave function
\begin{equation}\label{eq2}
\Psi(h;\mu_1,\mu_2)=\half \left[\Phi_1(-h;\mu_1)+\Phi_2(h;\mu_2)\right]: \H\times (\M_1\times\M_2)\to \bR
\end{equation}
has a saddle point ($\min$ in $h\in \H$, $\max$ in $(\mu_1,\mu_2)\in \M_1\times \M_2$).}
 \end{quote}
We associate with a saddle point $(h_*;\mu^*_1,\mu^*_2)$ the following entities:
\begin{itemize}
\item the {\sl risk}
\begin{equation}\label{eq11}
\epsilon_\star=\exp\{\Psi(h_*;\mu_1^*,\mu_2^*)\};
\end{equation}
 this quantity is uniquely defined by the saddle point value of $\Psi$ and thus is independent of how we select a saddle point;
\item the {\sl detector} $\phi_*(\omega)$ -- the affine function of $\omega\in\bR^d$ given by
\begin{equation}\label{eq3}
\phi_*(\omega)=h_*^T\omega+a,\,\,a=\half \left[\Phi_1(-h_*;\mu_1^*)-\Phi_2(h_*;\mu_2^*)\right].
\end{equation}
\end{itemize}
{A} simple sufficient condition for existence of a saddle point of (\ref{eq2}) is
\paragraph{Condition A:} {\sl The sets $\M_1$ and $\M_2$ are compact, and the function
$$
\overline{\Phi}(h)=\max_{\mu_1\in\M_1,\mu_2\in \M_2}\Phi(h;\mu_1,\mu_2)
$$
is coercive on $\H$, meaning that $\overline{\Phi}(h_i)\to\infty$ along every sequence $h_i\in\H$ with $\|h_i\|_2\to\infty$ as $i\to\infty$.}\\
\begin{quote}{\small
Indeed, under Condition A by Sion-Kakutani Theorem \cite{LH1993} it holds
$$
\SadVal[\Phi]:=\inf_{h\in \H}\underbrace{\max_{\mu_1\in M_1,\mu_2\in\M_2}\Phi(h;\mu_1,\mu_2)}_{\overline{\Phi}(h)}=\sup_{\mu_1\in M_1,\mu_2\in\M_2}\underbrace{\inf_{h\in \H}\Phi(h;\mu_1,\mu_2)}_{\underline{\Phi}(\mu_1,\mu_2)},
$$
 so that the optimization problems
$$
\begin{array}{ll}
(P):~~&{\Opt(P)=}\min\limits_{h\in \H}\overline{\Phi}(h)\\
(D):~~&{\Opt(D)=}\max\limits_{\mu_1\in\M_1,\mu_2\in\M_2}\underline{\Phi}(\mu_1,\mu_2)\\
\end{array}
$$
have equal optimal values. Under Condition A, problem $(P)$ clearly is a problem of minimizing a continuous coercive  function over a closed  set and as such is solvable; thus, $\Opt(P)=\Opt(D)$ is a real. Problem $(D)$ clearly is the problem of maximizing over a compact set of an upper semi-continuous  (since $\Phi$ is continuous) function taking real values and, perhaps, value $-\infty$, and not identically equal to $-\infty$ (since $\Opt(D)$ is a real), and thus $(D)$ is solvable. Thus, $(P)$ and $(D)$ are solvable with common optimal value, and therefore $\Phi$ has a saddle point.
}
\end{quote}
\subsection{Pairwise testing regular families of distributions}\label{sect:main}
\subsubsection{Main observation}
An immediate (and crucial!) observation is as follows:
\begin{proposition}\label{prop1}
In the situation of section \ref{sect:situ} and under Assumption I, one has
\begin{equation}\label{eq100}
\begin{array}{rcl}
\int_\Omega\exp\{-\phi_*(\omega)\}P(d\omega)&\leq&\epsilon_\star\,\,\forall P\in \P_1=\R[\H,\M_1,\Phi_1]\\
\int_\Omega\exp\{\phi_*(\omega)\}P(d\omega)&\leq&\epsilon_\star\,\,\forall P\in \P_2=\R[\H,\M_2,\Phi_2].\\
\end{array}
\end{equation}
\end{proposition}
{\bf Proof.} For every $\mu_1\in \M_1$, we have $\Phi_1(-h_*;\mu_1)\leq \Phi_1(-h_*;\mu_1^*)$, and for every $P\in\P_1$, we have
$$
\int_\Omega\exp\{-h_*^T\omega\}P(d\omega)\leq \exp\{\Phi_1(-h_*;\mu_1)\}
$$
for properly selected $\mu_1\in\M_1$. Thus,
$$
\int_\Omega\exp\{-h_*^T\omega\}P(d\omega)\leq \exp\{\Phi_1(-h_*;\mu_1^*)\}\,\,\forall P\in\P_1,
$$
whence also
$$
\int_\Omega\exp\{-h_*^T\omega-a\}P(d\omega)\leq \exp\{\half \left[\Phi_1(-h_*;\mu_1^*)+\Phi_2(h_*,\mu_2^*)\right]\}=\epsilon_\star\,\,\forall P\in\P_1.
$$
Similarly,
for every $\mu_2\in \M_2$, we have $\Phi_2(h_*;\mu_2)\leq \Phi_2(h_*;\mu_2^*)$, and for every $P\in\P_2$, we have
$$
\int_\Omega\exp\{h_*^T\omega\}P(d\omega)\leq \exp\{\Phi_2(h_*;\mu_2)\}
$$
for properly selected $\mu_2\in\M_2$. Thus,
$$
\int_\Omega\exp\{h_*^T\omega\}P(d\omega)\leq \exp\{\Phi_2(h_*;\mu_2^*)\}\,\,\forall P\in\P_2,
$$
whence
$$
\int_\Omega\exp\{h_*^T\omega+a\}P(d\omega)\leq \exp\{\half \left[\Phi_1(-h_*;\mu_1^*)+\Phi_2(h_*,\mu_2^*)\right]\}=\epsilon_\star\,\,\forall P\in\P_2.\eqno{\hbox{\qed}}
$$
{\subsubsection{Testing pairs of hypotheses}
\label{sect:pairwise}
\paragraph{Repeated observation.} Given $\Omega=\bR^d$, let random observations $\omega_t\in\Omega$, $t=1,2,...$, be generated as follows:
\begin{quote}
``In the nature'' there exists  a random sequence $\zeta_t\in\bR^N$, $t=1,2,...,$, of  {\sl driving factors} such that $\omega_t$ is a deterministic function of $\zeta^t=(\zeta_1,...,\zeta_t)$, $t=1,2,...$;
\end{quote}
we refer to this situation as to case of {\sl repeated observation $\omega^\infty$.}\par

Let now $(\H,\M,\Phi)$ be a regular data with observation space $\Omega$. We associate with this data four hypotheses
on the stochastic nature of observations $\omega^\infty=\{\omega_1,\omega_2,...\}$. Denoting  by $P_{|\zeta^{t-1}}$ the conditional, $\zeta^{t-1}$ being given, distribution of $\omega_t$, we say that the (distribution of the) repeated observation $\omega^\infty$ obeys hypothesis
\begin{itemize}
\item  $\HypR[\H,\M,\Phi]$, if $P_{|\zeta^{t-1}}\in\R[\H,\M,\Phi]$ for all $t$ and all $\zeta^{t-1}$;
\item $\HypS[\H,\M,\Phi]$, if $P_{|\zeta^{t-1}}\in\S[\H,\M,\Phi]$ for all $t$ and all $\zeta^{t-1}$;
\item $\HypRs[\H,\M,\Phi]$, if $P_{|\zeta^{t-1}}$ is independent of $t$ and $\zeta^{t-1}$ and belongs to $\R[\H,\M,\Phi]$;
\item $\HypSs[\H,\M,\Phi]$, if $P_{|\zeta^{t-1}}$ is independent of $t$ and $\zeta^{t-1}$ and belongs to $\S[\H,\M,\Phi]$.
\end{itemize}
Note that the last two hypotheses, in contrast to the first two, require from the observations $\omega_1,\omega_2,...$ to be i.i.d. Note also that $\HypR$ is weaker than $\HypS$, $\HypRs$ is weaker than $\HypSs$, and the ``non-stationary'' hypotheses  $\HypR$, $\HypS$ are weaker than their respective stationary counterparts $\HypRs$, $\HypSs$.
\par
The tests to be considered in the sequel operate with the initial fragment $\omega^K=(\omega_1,...,\omega_K)$, of a prescribed length $K$,  of the repeated observation $\omega^\infty=\{\omega_1,\omega_2,...\}$. We call $\omega^K$ {\sl $K$-repeated observation} and say that (the distribution of) $\omega^K$ obeys one of the above hypotheses, if $\omega^K$ is cut off the repeated observation $\omega^\infty$ distributed according to the hypothesis in question. We can think of $\HypR,\,\HypS,\,\HypRs$ and $\HypSs$ as hypotheses about the distribution of the $K$-repeated observation $\omega^K$.
\paragraph{Pairwise hypothesis testing from repeated observations.} Assume we are given two collections of regular data $(\H,\M_\chi,\Phi_\chi)$, $\chi=1,2$, with common observation space $\Omega=\bR^d$ and common $\H$. Given positive integer $K$ and $K$-repeated observation $\omega^K=(\omega_1,...,\omega_K)$, we want to decide on the pair of hypotheses $H_\chi=\HypR[\H,\M_\chi,\Phi_\chi]$, $\chi=1,2$, on the distribution of $\omega^K$.
\par
Assume that the convex-concave function (\ref{eq2}) associated with the pair of regular data in question has a saddle point $(h_*;\mu_1^*,\mu_2^*)$, and let $\phi_*(\cdot)$, $\epsilon_\star$ be the induced by this saddle point affine detector and its risk, see (\ref{eq3}), (\ref{eq11}). Let us set
$$
\phi_*^{(K)}(\omega^K)=\sum_{t=1}^K \phi_*(\omega_t).
$$
Consider decision rule $\T_*^K$ for hypotheses $H_\chi$, $\chi=1,2$, which, given observation $\omega^K$,
\begin{itemize}
\item accepts $H_1$ (and rejects $H_2$) if $\phi_*^{(K)}(\omega^K)>0$;
\item accepts $H_2$ (and rejects $H_1$) if $\phi_*^{(K)}(\omega^K)<0$;
\item in the case of a tie (when $\phi_*^{(K)}(\omega^K)=0$) the test, say, accepts $H_1$ and rejects $H_2$.
\end{itemize}
In what follows, we refer to $\T_*^K$ as to {\em test associated with detector  $\phi_*^{(K)}$.}
\begin{proposition}\label{prop2} In the situation in question, we have
\begin{equation}\label{mult12}
\begin{array}{lll}
(a)&\bE_{\zeta^K}\left\{\exp\{-\phi_*^{(K)}(\omega^{K})\}\right\}\leq\epsilon_\star^K,&\hbox{when $\omega^\infty$ obeys $H_1$};\\
(b)&\bE_{\zeta^K}\left\{\exp\{\phi_*^{(K)}(\omega^{K})\}\right\}\leq\epsilon_\star^K,&\hbox{when $\omega^\infty$ obeys $H_2$}.\\
\end{array}
\end{equation}
As a result, the test  $\T_*^K$ accepts exactly one of the hypotheses $H_1$, $H_2$, and the risk of this test --- the maximal, over $\chi=1,2$, probability not to accept the hypothesis $H_\chi$ when it is true (i.e., when the $K$-repeated observation $\omega^K$ obeys the  hypothesis $\HypR[\H,\M_\chi,\Phi_\chi]$) --- does not exceed $\epsilon_\star^K$.
\end{proposition}
{\bf Proof.} The fact that the test always accepts exactly one of the hypotheses $H_\chi$, $\chi=1,2$, is evident. Let us denote $\bE_{\zeta^t}$ the expectation w.r.t. the distribution of $\zeta^t$, and let $\bE_{\zeta_{t+1}|\zeta^t}$ stand for expectation w.r.t. conditional to $\zeta^t$ distribution of $\zeta_{t+1}$. Assuming that $H_1$ holds true and invoking the first inequality in (\ref{eq100}), we have
\bse
\bE_{\zeta_1}\{\exp\{-\phi_*^{(1)}(\omega_1)\}\}&\leq&\epsilon_\star,\\
\bE_{\zeta^{t+1}}\{\exp\{-\phi_*^{(t+1)}(\omega^{t+1})\}&=&\bE_{\zeta^{t}}\left\{\exp\{-\phi_*^{(t)}(\omega^{t})\}\bE_{\zeta_{t+1}|\zeta^t}
\{\exp\{-\phi_*(\omega_{t+1})\}\}\right\}\\
&\leq& \epsilon_\star\bE_{\zeta^{t}}\left\{\exp\{-\phi_*^{(t)}(\omega^{t})\}\right\}, \;1\leq t < K
\ese
(we have taken into account that $\omega_{t+1}$ is a deterministic function of $\zeta^{t+1}$ and that we are in the case where the conditional to $\zeta^t$ distribution of $\omega_{t+1}$ belongs to $\P_1=\R[\H,\M_1,\Phi_1]$),
and we arrive at (\ref{mult12}.$a$). This inequality
clearly implies that the probability to reject $H_1$ when the hypothesis is true is $\leq\epsilon_\star^K$ (since $\phi_*^{(K)}(\omega^K)\leq0$ when $\T_*^K$ rejects $H_1$).
 Assuming that $H_2$ is true and invoking the second inequality in (\ref{eq100}), similar reasoning shows that (\ref{mult12}.$b$)
 holds true, so that the probability to reject $H_2$ when the hypothesis is true does not exceed $\epsilon_\star^K$. \hfill\qed
\subsubsection{Illustration: sub-Gaussian and Gaussian cases}\label{subGaussian}
For $\chi=1,2$, let $U_\chi$ be nonempty closed convex set in $\bR^d$, and $\U_\chi$ be a compact convex subset of the interior of the positive semidefinite cone $\bS^d_+$. We assume that $U_1$ is compact. Setting
\begin{equation}\label{eq50}
\H_\chi=\Omega=\bR^d,\,\M_\chi=U_\chi\times\U_\chi,\;\;\Phi_\chi(h;\theta,\Theta)=\theta^Th+\half h^T\Theta h:\H_\chi\times \M_\chi\to\bR,\chi=1,2,
\end{equation}
we get two collections $(\H,\M_\chi,\Phi_\chi)$, $\chi=1,2$, of regular data. As we know from section \ref{Example1}, for $\chi=1,2$, the families of distributions  $\S[\bR^d,\M_\chi,\Phi_\chi]$ contain the families $\SG[U_\chi,\U_\chi]$ of sub-Gaussian distributions on $\bR^d$ with sub-Gaussianity parameters $(\theta,\Theta)\in U_\chi\times\U_\chi$ (see (\ref{eq40})), as well as families $\G[U_\chi,\U_\chi]$ of Gaussian distributions on $\bR^d$ with parameters $(\theta,\Theta)$ (expectation and covariance matrix) running through $U_\chi\times\U_\chi$. Besides this, the pair of regular data in question clearly satisfies
Condition A. Consequently, the test $\T_*^K$ given by the above construction as applied to the  collections of regular data (\ref{eq50}) is well defined and allows to decide on hypotheses $H_\chi=\HypR[\bR^d,U_\chi,\U\chi]$, $\chi=1,2$, on the distribution of the $K$-repeated observation $\omega^K$. The same test can be also used to decide on stricter hypotheses $H^G_\chi$, $\chi=1,2$, stating that the observations $\omega_1,...,\omega_K$ are i.i.d. and drawn from a Gaussian distribution $P$ belonging to $\G[U_\chi,\U_\chi]$. Our goal now is to process in  detail the situation in question and to refine our conclusions on the risk of the test $\T_*^1$ when the {\em Gaussian} hypotheses $H^G_\chi$ are considered and the situation is {\sl symmetric}, {that is,} when $\U_1=\U_2$.
\par
Observe, first, that the convex-concave function $\Psi$ from (\ref{eq2}) in the situation under consideration becomes
\begin{equation}\label{Psi}
\Psi(h;\theta_1,\Theta_1,\theta_2,\Theta_2)=\half h^T[\theta_2-\theta_1]+{1\over 4}h^T\Theta_1h + {1\over 4} h^T\Theta_2h.
\end{equation}
We are interested in solutions to the saddle point problem -- find a saddle point of function (\ref{Psi}) --
 \begin{equation}\label{PsiP}
 \min_{h\in\bR^d}\max_{\stackrel{\theta_1\in U_1,\theta_2\in U_2}{\Theta_1\in\U_1,\Theta_2\in\U_2}}\Psi(h;\theta_1,\Theta_1,\theta_2,\Theta_2)
 \end{equation}
 From the structure of $\Psi$ and compactness of $U_1$, $\U_1$, $\U_2$, combined with the fact that $\U_\chi$, $\chi=1,2$, are comprised of positive definite matrices, it immediately follows that saddle points do exist, and a saddle point $(h_*;\theta_1^*,\Theta_1^*,\theta_2^*,\Theta_2^*)$ satisfies the relations
\begin{equation}\label{eq41}
\begin{array}{ll}
(a)&h_*=[\Theta_1^*+\Theta_2^*]^{-1}[\theta_1^*-\theta_2^*],\\
(b)&h_*^T(\theta_1-\theta_1^*)\geq0\;\forall \theta_1\in U_1, \;\;h_*^T(\theta_2^*-\theta_2)\geq 0\;\forall \theta_2\in U_2,\\
(c)& h_*^T\Theta_1h_*\leq h_*^T\Theta_1^*h_*\;\forall \Theta_1\in\U_1,\;\;h_*^T\Theta_2^*h_*\leq h_*\Theta_2^*h_*\;\forall \Theta_2\in\U_2.
\end{array}
\end{equation}
From (\ref{eq41}.$a$) it immediately follows that the affine detector $\phi_*(\cdot)$ and risk $\epsilon_\star$, as given by \rf{eq11} and \rf{eq3}, are
\begin{equation}\label{eq333}
\begin{array}{rcl}
\phi_*(\omega)&=&h_*^T[\omega-w_*]+\half h_*^T[\Theta_1^*-\Theta_2^*]h_*,\,\,
w_*=\half [\theta_1^*+\theta_2^*];\\
\epsilon_\star&=&\exp\{-{1\over 4}[\theta_1^*-\theta_2^*]^T[\Theta_1^*+\Theta_2^*]^{-1}[\theta_1^*-\theta_2^*]\}\\
&=&\exp\{-{1\over 4}h_*^T[\Theta_1^*+\Theta_2^*]h_*\}.\\
\end{array}
\end{equation}
Note that in the {\sl symmetric case} (where $\U_1=\U_2$), there always exists a saddle point of $\Psi$ with $\Theta_1^*=\Theta_2^*$, and the test $\T_*^1$ associated with such  saddle point is quite transparent: it is the maximum likelihood test for two Gaussian distributions, $\N(\theta_1^*,\Theta_*)$,
$\N(\theta_2^*,\Theta_*)$, where $\Theta_*$ is the common value of $\Theta_1^*$ and $\Theta_2^*$, and the bound $\epsilon_\star$ for the risk of the test is nothing but the Hellinger affinity of these two Gaussian distributions, or, equivalently,
\begin{equation}\label{epsilonstarsubG}
\epsilon_\star=\exp\left\{-{\eight}[\theta_1^*-\theta_2^*]^T\Theta_*^{-1}[\theta_1^*-\theta_2^*]\right\}.
\end{equation}
Applying Proposition \ref{prop2}, we arrive at the following result:
\begin{proposition}\label{propsubGauss}
In the symmetric sub-Gaussian case (i.e., in the case of {\rm (\ref{eq50})} with $\U_1=\U_2$), saddle point problem {\rm (\ref{Psi}), (\ref{PsiP})} admits a saddle point of the form
$(h_*;\theta_1^*,\Theta_*,\theta_2^*,\Theta_*)$, and the associated affine detector and its risk are given by
\begin{equation}\label{eq3333}
\begin{array}{rcl}
\phi_*(\omega)&=&h_*^T[\omega-w_*],\,\,
w_*=\half [\theta_1^*+\theta_2^*];\\
\epsilon_\star&=&\exp\{-{1\over 8}[\theta_1^*-\theta_2^*]^T\Theta_*^{-1}[\theta_1^*-\theta_2^*]\}.\\
\end{array}
\end{equation}
As a result, when deciding, via $\omega^K$, on ``sub-Gaussian hypotheses'' $\HypS[\bR^d,\M_\chi,\M_\chi]$, $\chi=1,2$ (in fact - even on weaker hypotheses $\HypR[\bR^d,\M_\chi,\M_\chi]$, $\chi=1,2$), the risk of the test $\T_*^K$ {associated with $\phi_*^{(K)}(\omega^K):=\sum_{t=1}^K\phi_*(\omega_t)$} is at most $\epsilon_\star^K$.
\end{proposition}
 In the symmetric single-observation Gaussian case, that is, when $\U_1=\U_2$ and we apply the test $\T_*=\T_*^1$ to observation $\omega\equiv\omega_1$ in order to decide on the hypotheses $\H^G_\chi$, $\chi=1,2$, the above risk bound can be improved:
\begin{proposition}\label{prop22}
Consider symmetric case $\U_1=\U_2=\U$, let $(h_*;\theta_1^*;\Theta_1^*,\theta_2^*,\Theta_2^*)$ be ``symmetric'' -- with $\Theta_1^*=\Theta_2^*=\Theta_*$ --
 saddle point of function $\Psi$ given by {\rm (\ref{Psi})}, and let $\phi_*$ be the affine detector given by {\rm (\ref{eq41}) and (\ref{eq333})}:
$$
\phi_*(\omega)=h_*^T[\omega-w_*],\;\;h_*={1\over2}\Theta_*^{-1}[\theta_1^*-\theta_2^*],\,\,w_*=\half [\theta_1^*+\theta_2^*].
$$
Let also
\begin{equation}\label{delta0}
\delta=\sqrt{h_*^T\Theta_*h_*}=\half\sqrt{[\theta_1^*-\theta_2^*]^T\Theta_*^{-1}[\theta_1^*-\theta_2^*]},
\end{equation}
so that
\begin{equation}\label{delta}
\delta^2=h_*^T[\theta_1^*-w_*]=h_*^T[w_*-\theta_2^*]\;\mbox{\rm and}\; \epsilon_\star=\exp\{-\half \delta^2\}.
\end{equation}
 Let, further, $\alpha\leq\delta^2$, $\beta\leq\delta^2$.
Then
\begin{equation}\label{eq1100}
\begin{array}{ll}
(a)&\forall (\theta\in U_1,\Theta\in \U): \Prob_{\omega\sim \N(\theta,\Theta)}\{\phi_*(\omega)\leq\alpha\}\leq \Erf(\delta-\alpha/\delta)\\
(b)&\forall (\theta\in U_2,\Theta\in \U): \Prob_{\omega\sim \N(\theta,\Theta)}\{\phi_*(\omega)\geq-\beta\}\leq \Erf(\delta-\beta/\delta),\\
\end{array}
\end{equation}
where
$$
\Erf(s)={1\over\sqrt{2\pi}}\int_s^\infty\exp\{-r^2/2\}dr
$$
is the normal error function. In particular, when deciding, via a single observation $\omega$, on Gaussian hypotheses $H^G_\chi$, $\chi=1,2$, with $H_\chi^G$ stating that $\omega\sim\N(\theta,\Theta)$ with $(\theta,\Theta)\in U_\chi\times\U$,
the risk of the test {$\T_*^1$} {associated with $\phi_*$} is at most $\Erf(\delta)$.
\end{proposition}
{\bf Proof.}
Let us prove $(a)$ (the proof of $(b)$ is completely similar). For $\theta\in U_1$, $\Theta\in \U$ we have
$$
\begin{array}{l}
\Prob_{\omega\sim \N(\theta,\Theta)}\{\phi_*(\omega)\leq\alpha\}=\Prob_{\omega\sim \N(\theta,\Theta)}\{h_*^T[\omega-w_*]\leq\alpha\}\\
=\Prob_{\xi\sim\N(0,I)}\{h_*^T[\theta+\Theta^{1/2}\xi-w_*]\leq\alpha\}\\
=\Prob_{\xi\sim\N(0,I)}\{[\Theta^{1/2}h_*]^T\xi\leq\alpha-\underbrace{h_*^T[\theta-w_*]}_{{\geq h_*^T[\theta_1^*-w_*]=\delta^2\atop
\hbox{\ \tiny by (\ref{eq41}.$b$),(\ref{delta})}}}\}
\leq \Prob_{\xi\sim \N(0,I)}\{[\Theta^{1/2}h_*]^T\xi\leq\alpha-\delta^2\}\\
=\Erf([\delta^2-\alpha]/\|\Theta^{1/2}h_*\|_2)\\
\leq \Erf([\delta^2-\alpha]/\|\Theta_*^{1/2}h_*\|_2)\hbox{\ [since $\delta^2-\alpha\geq0$ and $h_*^T\Theta h_*\leq h_*^T\Theta_*h_*$ by (\ref{eq41}.$c$)]}\\
=\Erf([\delta^2-\alpha]/\delta).\\
\end{array}
$$
The ``in particular'' part of Proposition is readily given by (\ref{eq1100}) as applied with $\alpha=\beta=0$. \hfill $\Box$
}
\subsection{Testing multiple hypotheses from repeated observations}\label{sect:multiple}
Consider  the situation as follows: we are given
\begin{itemize}
\item observation space $\Omega=\bR^d$ and a symmetric w.r.t. the origin closed convex set $\H$;
\item $J$ closed convex sets $\M_j\subset\bR^{n_j}$, $j=1,...,J$, along with $J$ convex-concave continuous functions $\Phi_j(h;\mu^j):\H\times\M_j\to\bR$.
\end{itemize}
These data give rise to $J$ hypotheses $H_j=\HypR[\H,\M_j,\Phi_j]$ on the distribution of repeated observation $\omega^\infty=\{\omega_t\in\bR^d,t\geq1\}$. On the top of it, assume we are given a {\sl closeness} relation -- a subset $\C\subset\{1,...,J\}^2$ which we assume to contain the diagonal ($(j,j)\in \C$ for all $j\leq J$) and to be symmetric ($(i,j)\in\C$ if and only if $(j,i)\in\C$). In the sequel, we interpret indexes $i,j$ with $(i,j)\in\C$ (and the hypotheses $H_i$, $H_j$ with these indexes) as {\sl $\C$-close to each other}.
\par
Our goal is,  given a positive integer $K$ and $K$-repeated observation $\omega^K=(\omega_1,...,\omega_K)$, to decide,  ``up to closeness $\C$,'' on the hypotheses $H_j$, $j=1,...,J$ (which is convenient to be thought of as hypotheses about the distribution of $\omega^K$).
\par
Let us act as follows\footnote{The construction we are about to present and the first related result (Proposition \ref{prop3}) originate from \cite[Section 3]{GJN}; we reproduce them below to make our exposition self-contained.}. {Let us make}
 \begin{quote}
 {\bf Assumption II} {\sl For every pair $i,j$, $1\leq i<j\leq J$, with $(i,j)\not\in\C$, the convex-concave function
$$
\Psi_{ij}(h;\mu_i,\mu_j)=\half \left[\Phi_i(-h;\mu_i)+\Phi_j(h;\mu_j)\right]: \;\H\times (\M_i\times\M_j)\to \bR
$$
has a saddle point $(h^{ij};\mu^{ij}_i,\mu^{ij}_j)$ on $\H\times (\M_i\times\M_j)$ ($\min$ in $h$, $\max$ in $\mu_i,\mu_j$)},
\end{quote}
 and let us set
\begin{equation}\label{assume!}
\left.\begin{array}{rcl}
\epsilon_{ij}&=&\exp\{\Psi_{ij}(h^{ij};\mu_i^{ij},\mu_j^{ij})\},\\
\phi_{ij}(\omega)&=&[h^{ij}]^T\omega+a_{ij}\\
&&a_{ij}=\half \left[\Phi_i(-h^{ij};\mu_i^{ij})-\Phi_j(h^{ij};\mu_j^{ij})\right]\\
\end{array}\right\}, 1\leq i<j\leq J\ \&\ (i,j)\not\in\C
\end{equation}
(cf. (\ref{eq11}) and (\ref{eq3})).
We further set
$$
\phi_{ij}(\cdot)\equiv 0, \epsilon_{ij}=1 \,\,\forall (i,j)\in\C,
$$
and
$$
\phi_{ij}(\omega)=-\phi_{ji}(\omega),\,\epsilon_{ij}=\epsilon_{ji},\,1\leq j < i\leq J\ \&\ (i,j)\not\in\C,
$$
and set
$$
\phi_{ij}^{(K)}(\omega^K)=\sum_{t=1}^K\phi_{ij}(\omega_t),\,1\leq i,j\leq J.
$$
Now, by construction,
\[
\phi_{ij}^{(K)}(\cdot)=-\phi_{ji}^{(K)}(\cdot),\,\epsilon_{ij}=\epsilon_{ji},\,1\leq i,j\leq J.\\
\]
Observe that for every $j\leq J$ such that the $K$-repeated observation $\omega^K$ obeys hypothesis $H_j=\HypR[\H,\M_j,\Phi_j]$, we have
\begin{equation}\label{eqrisk}
\bE\left\{\exp\{\phi_{ij}^{(K)}(\omega^K)\}\right\}\leq\epsilon_{ij}^K,\,i=1,...,J.
\end{equation}
(we have used (\ref{mult12}) along with $\phi_{ij}\equiv -\phi_{ji}$).
\par
Furthermore, whenever $[\alpha_{ij}]_{i,j\leq J}$ is a skew-symmetric matrix (i.e., $\alpha_{ij}=-\alpha_{ji}$), the shifted detectors
\begin{equation}\label{eq20}
\widehat{\phi}_{ij}^{(K)}(\omega^K)=\phi_{ij}^{(K)}(\omega^K)+\alpha_{ij}\\
\end{equation}
satisfy the scaled version of \rf{eqrisk}, specifically,
for every $j\leq J$ such that the $K$-repeated observation $\omega^K$ obeys hypothesis $H_j=\HypR[\bR^d,\M_j,\Phi_j]$, we have
\begin{equation}\label{eqrisk1}
\bE\left\{\exp\{\widehat{\phi}_{ij}^{(K)}(\omega^K)\}\right\}\leq\epsilon_{ij}^K\exp\{\alpha_{ij}\},\,i=1,2,...,J.
\end{equation}
The bottom line is as follows:
\begin{proposition}\label{prop3} In the situation in question, given a closeness $\C$, consider the following test $\T_{\C}^K$ deciding on the hypotheses $H_j$, $1\leq j\leq J$, on the distribution of $K$-repeated observation $\omega^K$: given skew-symmetric shift matrix $[\alpha_{ij}]_{i,j}$ and observation $\omega^K$, $\T_{\C}^K$ accepts all hypotheses $H_i$ such that
$$
\widehat{\phi}^{(K)}_{ij}(\omega^K)>0 \hbox{\ whenever\ } (i,j)\not\in\C\eqno{(*_i)}
$$
and reject all hypotheses $H_i$ for which the predicate $(*_i)$ does not take place. Then
\par
{\rm (i)} Test $\T_{\C}^K$
 accepts some of (perhaps, none of) hypotheses $H_i$, $i=1,...,J$, and all accepted hypotheses, if any, are $\C$-close to each other. Besides, $\T_{\C}^K$ has $\C$-risk at most
$$
\widehat{\epsilon}=\max\limits_i\sum\limits_{j:(i,j)\not\in\C}\epsilon_{ij}^K\exp\{-\alpha_{ij}\},
$$
meaning that for every $j_*\leq J$ such that the distribution ${\bar{P}_K}$ of  $\omega^K$ obeys the hypothesis $H_{j_*}$ (i.e., the hypothesis $H_{j_*}$ is true), the ${\bar{P}_K}$-probability of the event
\begin{quote}
''either the true hypothesis $H_{j_*}^K$ is not accepted, or the list of accepted hypotheses contains a  hypothesis which is not $\C$-close to $H_{j_*}$''
\end{quote}
does not exceed $\widehat{\epsilon}$.
\par
{\rm (ii)} The infimum of $\widehat{\epsilon}$ over all skew-symmetric shifts $\alpha_{ij}$ is exactly the spectral norm $
\|E^{(K)}\|_{2,2}$ of the symmetric entry-wise nonnegative matrix
\be
E^{(K)}=\left[E_{ij}^{(K)}=\left\{\begin{array}{ll}\epsilon_{ij}^K,&(i,j)\not\in\C\\
0,&(i,j)\in\C\\
\end{array}\right.\right]_{i,j=1}^J.
\ee{E^Kmatrix}
This infimum is attained when the Perron-Frobenius eigenvector $g$ of $E^{(K)}$ can be selected to be positive, in which case an optimal selection of $\alpha_{ij}$ is the selection
$$
\alpha_{ij}=\ln(g_i/g_j),\;1\leq i,j\leq J.
$$
\end{proposition}
{\bf Proof.} Given $\omega^K$, the test $\T_{\C}^K$ can accept hypotheses $H_i$ and $H_j$ with $(i,j)\not\in\C$ only when $\phi^{(K)}_{ij}(\omega^K)>0$ and $\phi^{(K)}_{ji}(\omega^K)>0$, which is impossible due to $\phi_{ij}^{(K)}(\cdot)=-\phi_{ji}^{(K)}(\cdot)$. Thus, $\T_{\C}^K$ can accept $H_i$ and $H_j$ only when $(i,j)\in\C$. Further, let the distribution ${\bar{P}_K}$ of $\omega^K$ obey hypothesis $H_{j_*}$. Invoking (\ref{eqrisk1}) and the relation {
$\widehat{\phi}_{j_*j}^{(K)}(\cdot)= -\widehat{\phi}_{jj_*}^{(K)}(\cdot)$}, for every $j\leq J$ with $(j_*,j)\not\in\C$,
the ${\bar{P}_K}$-probability of the event $\widehat{\phi}_{j_*j}^{(K)}(\omega^K)\leq 0$, or, which is the same, of the event  $\widehat{\phi}_{jj_*}^{(K)}(\omega^K)\geq 0$, is at most $\epsilon_{jj_*}^K\exp\{\alpha_{jj_*}\}=\epsilon_{j_*j}^K\exp\{-\alpha_{j_*j}\}$. Using the union bound, the ${\bar{P}_K}$-probability of the event ``$H_{j_*}$ is not accepted'' is at most
$$\sum_{j:(j_*j)\not\in \C} \epsilon_{j_*j}^K\exp\{-\alpha_{j_*j}\}=\sum_{j}E^{(K)}_{j_*j}\exp\{-\alpha_{j_*j}\}\leq\widehat{\epsilon}.
 $$
 By construction of the test, when $H_{j_*}$ is accepted and $j$ is not $\C$-close to $j_*$, $H_j$ is not accepted {(we have already seen that the test never accepts a pair of hypotheses which are not $\C$-close to each other)}. Thus, the $\C$-risk of $\T^K_{\C}$ indeed is at most $\widehat{\epsilon}$. Now, $E^{(K)}$ is a symmetric entry-wise nonnegative matrix, so that its leading eigenvector $g$ can be selected to be nonnegative. When $g$ is positive, setting
$
\alpha_{ij} =\ln(g_i/g_j),
$
we get for every $i$
$$
\sum_{j}E^{(K)}_{ij}\exp\{-\alpha_{ij}\}=\sum_{j}E^{(K)}_{ij}g_j/g_i=(E^{(K)}g)_i/g_i=\|E^{(K)}\|_{2,2},
$$
and thus $\widehat{\epsilon}=\|E^{(K)}\|_{2,2}$. The fact that this is the smallest possible, over skew-symmetric shifts $\alpha_{ij}$, value of $\widehat{\epsilon}$ is proved in \cite{GJN}. When $g$ is just nonnegative, consider close to $E^{(K)}$ symmetric matrix $\widehat{E}$ with positive entries $\widehat{E}_{ij}\geq E^{(K)}_{ij}$; utilizing the (automatically strictly positive) Perron-Frobenius
eigenvector $g$ of this matrix, we, as was just explained, get skew-symmetric shifts $\alpha_{ij}$ such that
$$
\sum_j\widehat{E}_{ij}\exp\{-\alpha_{ij}\}\leq \|\widehat{E}\|_{2,2}
$$
for all $i$; the left hand side in this inequality is $\geq \sum_jE_{ij}^{(K)}\exp\{-\alpha_{ij}\}$, and the right hand side can be made arbitrarily close to $\|E^{(K)}\|_{2,2}$ by making $\widehat{E}$ close enough to $E^{(K)}$. Thus, we indeed can make $\widehat{\epsilon}$ arbitrarily close to $
\|E^{(K)}\|_{2,2}$.\hfill \qed
\subsubsection{Special case: inferring colors}\label{sect:IC}
 Assume that we are given $J$ collections $(\H,\M_j,\Phi_j)$, $1\leq j\leq J$, of regular data with common observation space $\Omega=\bR^d$ and common $\H$, and thus have at our disposal $J$ hypotheses $H_j=\HypR[\H,\M_j,\Phi_j]$, on the distribution of $K$-repeated observation $\omega^K$. Let also the index set $\{1,...,J\}$ be partitioned into $L\geq2$ non-overlapping nonempty subsets $\I_1,...,\I_L$; it is convenient to think about $\ell$, $1\leq \ell\leq L$, as the common color of indexes $j\in \I_\ell$, and  that the colors of indexes $j$ are inherited by the hypotheses $H_j$. The Color Inference ({CI}) problem we want to solve amounts to decide,
 given $K$-repeated observation $\omega^K$ obeying one or more  of the hypotheses $H_1,...,H_J$, on the color of these hypotheses. Note that it may happen that the distribution of $\omega^K$ obeys a pair of hypotheses $H_i$, $H_j$ of different colors. If it is not the case -- that is, no distribution of $\omega^K$ obeys a pair of hypotheses $H_i$, $H_j$ of two distinct  colors -- we call the CI problem well-posed.
 In the well-posed case, we can speak about the color of the distribution of $\omega^K$ provided this distribution obeys the union, over $j=1,...,J$, of hypotheses $H_j$; this is the color of (any of) the hypotheses $H_j$ obeyed by the distribution of $\omega^K$, and the CI problem is to infer this color given $\omega^K$.
 \par
In order to process the CI problem via our machinery, let us define closeness $\C$ as follows:
$$
(i,j)\in\C \Leftrightarrow \hbox{\ $i$ and $j$ are of the same color}.
$$
Assuming that the resulting $\C$ ensures validity of Assumption II, we can apply the above scheme to build test $\T_{\C}^K$. We can then convert this test into a color inference as follows. Given a $K$-repeated observation $\omega^K$,
it may happen that $\T_{\C}^K$, as applied to $\omega^K$, accepts one or more among the hypotheses $H_j$. In this case, by item (i) of Proposition \ref{prop3}, all accepted hypotheses are $\C$-close to each other (in other words, are of the same color), and we claim that this is the color of the distribution of $K$-repeated observation we are dealing with. And if $\T_{\C}^K$, as applied to $\omega^K$, accepts nothing, we claim that the color  we are interested in remains undecided.
\par
Let us analyze the just described color inferring procedure, let it be denoted $\A^K$.  Observe, first, that in the situation in question, assuming w.l.o.g. that the sets $\I_1,\I_2,...,\I_L$ are consecutive fragments in $\{1,...,J\}$, the matrix $E^{(K)}$ given by
 \rf{E^Kmatrix} is naturally partitioned into $L\times L$ blocks $E^{pq}=(E^{qp})^T$, $1\leq p,q\leq L$, where $E^{pq}$ is comprised of entries $E^{(K)}_{ij}$ with $i\in \I_p$, $j\in\I_q$. By construction of $E^{(K)}$, the diagonal blocks $E^{pp}$ are zero, and off-diagonal blocks are entry-wise positive (since $\epsilon_{ij}$ clearly is positive for all pairs $i,j$ of different colors).
 It follows that Perron-Frobenius eigenvectors of $E^{(K)}$ are strictly positive. This implies that for properly selected shifts $\alpha_{ij}=-\alpha_{ji}$, the quantity $\widehat{\epsilon}$ in Proposition \ref{prop3} is equal to $\|E^{(K)}\|_{2,2}$; in what follows we assume that the test $\T_{\C}^K$ utilizes exactly these optimal shifts, so that we are in the case of $\widehat{\epsilon}=\|E^{(K)}\|_{2,2}$.
  \par
  Now, it may happen (``bad case'') that
 that $\|E^{(K)}\|_{2,2}\geq 1$; in this case Proposition \ref{prop3} says nothing meaningful about the quality of the test $\T_{\C}^K$, and consequently, we cannot say much about the quality of $\A^K$. In contrast to this, we claim that
 \begin{lemma}\label{lem:diese} Assume that $\widehat{\epsilon}:=\|E^{(K)}\|_{2,2}<1$. Then the CI problem is well posed, and whenever
 the distribution  $\bar{P}_K$ of $\omega^K$ obeys one of the hypotheses $H_j$, $j=1,...,J$, $\A^K$ recovers correctly the color of $\bar{P}_K$ with $\bar{P}_K$-probability at least $1-\widehat{\epsilon}$.
 \end{lemma}
 \begin{quote}{\small
 The proof is immediate. In the good case, all entries in $E^{(K)}$ are of magnitude $<1$, whence $\epsilon_{ij}<1$ whenever $(i,j)\not\in\C$, see (\ref{E^Kmatrix}), so that
 \begin{equation}\label{goodcase}
 \bar{\epsilon}:=\max_{i,j}\{\epsilon_{ij}:(i,j)\not\in\C\}<1.
 \end{equation}
 It follows that the nonzero entries in $E^{(M)}$ are nonnegative and $\leq \bar{\epsilon}^M$, whence
 \begin{equation}\label{normbound}
 \widehat{\epsilon}(M):=\|E^{(M)}\|_{2,2}\leq J\bar{\epsilon}^M\to0,\,\mbox{ as $M\to\infty$}.
 \end{equation}
 In particular, for properly selected $M$ we have $\widehat{\epsilon}(M)<1/2$. Applying Proposition \ref{prop3} with $M$ in the role of $K$, we see that if the distribution ${\bar{P}_K}$ of $\omega^K$ obeys hypothesis $H_{j_*}$  with some $j_*\leq J$ (due the origin of our hypotheses, this is exactly the same as to say that the distribution ${\bar{P}_M}$ of $\omega^M$ obeys $H_{j_*}$), then with ${\bar{P}_M}$-probability at least $1-\widehat{\epsilon}(M)>1/2$ the test $\T_{\C}^M$ accepts hypothesis $H_{j_*}$.
 It follows that if ${\bar{P}_M}$ obeys both $H_{j'}$ and $H_{j''}$, then $\T_{\C}^M$ will accept $H_{j'}$ and $H_{j''}$ simultaneously with  positive ${\bar{P}_M}$-probability, and since, as we have already explained, $\T_{\C}^M$ never accepts two hypotheses of different color simultaneously, we conclude that $H_{j'}$ and $H_{j''}$ are of the same color. This conclusion holds true whenever the distribution of $\omega^K$ obeys one or more  of the hypotheses $H_j$, $1\leq j\leq K$, meaning that the CI  problem is well posed.
 \par
 Invoking Proposition \ref{prop3}, we conclude that if the distribution $\bar{P}_K$ of $\omega^K$ obeys, for some $j_*$, the hypothesis $H_{j_*}$, then the ${\bar{P}_K}$-probability for $\T_{\C}^K$ to accept $H_{j_*}$ is at least $1-\widehat{\epsilon}(K)$; and from the preceding analysis, whenever $\T_{\C}^K$ accepts $H_{j_*}$ such that ${\bar{P}_K}$ obeys $H_{j_*}$, $\A^K$ correctly infers the color of ${\bar{P}_K}$, as claimed. \hfill $\Box$
 }\end{quote}
 \par
 Finally, we remark that when (\ref{goodcase}) holds, (\ref{normbound}) implies that $\widehat{\epsilon}(K)\to 0$ as $K\to\infty$, so that the CI problem is well posed, and for every desired risk level $\epsilon\in(0,1)$ we can find efficiently observation time $K=K(\epsilon)$ such that $\widehat{\epsilon}(K)\leq\epsilon$. As a result, for this $K$ the color inferring procedure $\A_K$ recovers the color of the distribution ${\bar{P}_K}$ of $\omega^K$ (provided this distribution obeys some of the hypotheses $H_1,...,H_J$) with ${\bar{P}_K}$-probability at least $1-\epsilon$.

\section{Application: aggregating estimates by testing}\label{sect:aggr}
Let us consider the situation as follows:
\begin{itemize}
\item We are given $I$ triples of regular data $(\H,\M_i,\Phi_i)$, $1\leq i\leq I$, with common $\H$ and $\Omega=\bR^d$ and the parameter sets $\M_i$ sharing the common embedding space $\bR^n$; for the sake of simplicity, assume that $\M_i$ are bounded (and thus are nonempty convex compact sets in $\bR^n$) and the continuous convex-concave functions $\Phi_i(h;\mu):\H\times \M_j\to\bR$ are coercive in $H$: $\Phi_i(h_t,\mu)\to\infty$ whenever $\mu\in\M_i$ and sequence $\{h_t\in \H,t\geq 1\}$ satisfies $\|h_t\|_2\to\infty$, $t\to\infty$.
\item We observe a realization of $K$-repeated observation $\omega^K=(\omega_1,...,\omega_K)$ with i.i.d. $\omega_t$'s drawn from unknown probability distribution $\bar{P}$ known to belong to the family $\P=\bigcup\limits_{i\leq I}\S[\H,\M_i,\Phi_i]$. Thus, ``in the nature'' there exists $\bar{i}\leq I$ and $\bar{\mu}\in \M_{\bar{i}}$ such that
    \begin{equation}\label{eq1939}
    \ln\left(\bE_{\omega\sim \bar{P}}\{\exp\{h^T\omega\}\}\right)\leq \Phi_{\bar{i}}(h,\bar{\mu})\,\,\forall h\in\H.
    \end{equation}
    we call $\bar{\mu}$ the  parameter associated with $\bar{P}$, and our goal  is to recover from $K$-repeated observation $\omega^K$ the image $\bar{g}=G\bar{\mu}$ of $\bar{\mu}$ under a given linear mapping $\mu\mapsto G\mu:\bR^n\to\bR^m$.
\end{itemize}
Undoubtedly, parameter estimation problem is a fundamental problem of mathematical statistics, and as such is the subject of a huge literature. In particular, several constructions of estimators based on testing of convex hypotheses have been studied in connection with signal reconstruction \cite{Birge1983,Birge2006,Birge2013} and linear functionals estimation \cite{Donoho1987,donoho1991}.
Our actual goal to be addressed below is more modest:  we assume that we are given $L$ candidate estimates $g_1,...,g_L$ of $\bar{g}$ (these estimates could {be outputs of various estimation routines} applied to independent observations sampled form $\bar{P}$), and our goal is to select the best -- the ${\|\cdot\|_2}$-closest to $\bar{g}$ among the estimates $g_1,...,g_L$. This is the well-known problem of {\sl aggregating estimates}, and our goal is to process this aggregation problem  via the {Color Inference}
    procedure from section \ref{sect:IC}.
\subsection{Aggregation procedure}
 It should be stressed that as stated, the aggregation problem appears to be ill-posed: there could be several pairs $(\bar{i},\bar{\mu})$ satisfying (\ref{eq1939}), and the values of $G\mu$ at the $\mu$-components of these pairs could be different for different pairs, so that $\bar{g}$ not necessary is well defined. One way to resolve this ambiguity would be to assume that given $\bar{P}\in\P$, relation (\ref{eq1939}) uniquely defines $\bar{\mu}$. We, however, prefer another setting: $\bar{\mu}$ and $\bar{i}$ satisfying (\ref{eq1939}), same as $\bar{P}\in\P$, are ``selected by nature'' (perhaps, from several alternatives),  and the performance of the aggregating procedure we are about to develop should be independent of what is the nature's selection.
    \par
When processing the  aggregation problem, we assume w.l.o.g. that all points $g_1,...,g_L$ are distinct from each other.
Let us split the space $\bR^m$ where $G\mu$ takes values into $L$ {\sl Voronoi cells}
\begin{equation}\label{Voronoi}
\begin{array}{rcl}
V_\ell&=&\{g\in\bR^m: \,\|g-g_\ell\|_2\leq \|g-g_{\ell'}\|_2\;\forall \ell'\leq L\}\\
&=&\{g\in \bR^m:\, u_{\ell\ell'}^Tg\leq v_{\ell\ell'}\;\forall (\ell'\leq L,\ell'\neq\ell)\},\\
&&u_{\ell\ell'}=\|g_{\ell'}-g_\ell\|_2^{-1}[g_{\ell'}-g_\ell],\;\;v_{\ell\ell'}=\half u_{\ell\ell'}^T[g_{\ell'}+g_\ell],\;1\leq \ell,\ell'\leq L,\ell\neq\ell'.
\end{array}
\end{equation}
Note that $V_\ell$ is comprised of all points $g$ from $\bR^m$ for which $g_\ell$ is (one of) the $\|\cdot\|_2$-closest to $g$ among the points $g_1,...,g_L$.
Let us set
\begin{equation}\label{VoronoiW}
W^i_\ell=\{\mu\in\M_i:G\mu\in V_\ell\},\,1\leq i\leq I,1\leq\ell\leq L,
\end{equation}
so that $W^i_\ell$ are convex compact sets in $\bR^m$. Observe that $g_\ell$ can be a solution to the aggregation problem (that is, the closest to $\bar{g}$ point among $g_1,...,g_L$) only when $\bar{g}=G\bar{\mu}$ belongs to $V_\ell$, that is, only when $\bar{\mu}\in W^{\bar{i}}_\ell$ for some $\bar{i}$, implying that at least one of the sets $W^1_\ell,\,W^2_\ell,\,...,\,W^I_\ell$  is nonempty. Whether the latter condition is indeed satisfied for a given $\ell$ this can be found out efficiently via solving $I$ convex feasibility problems.
If the latter condition does {\sl not} hold for some $\ell$ (let us call the associated estimate $g_\ell$  {\sl redundant}), we can eliminate $g_\ell$ from the list of estimates to be aggregated without affecting the solution to the aggregation problem.
Then we can redefine the Voronoi cells for the reduced list of estimates in the role of our original list, check whether this list still contains a redundant estimate, eliminate the latter, if it exists, and proceed recursively until a list of estimates (which by construction still contains all solutions to the aggregation problem) with no redundant estimates is built. We lose nothing when assuming that this ``purification''  was carried out in advance, so that already the initial list $g_1,...,g_L$ of estimates does not contain redundant ones. Of course, we lose nothing when assuming that $L\geq 2$. Thus, from now on
we assume that $L\geq2$ and for every $\ell\leq L$, at least one of the sets $W^i_\ell$, $1\leq i\leq I$, is nonempty.
\par
Note that to solve the aggregation problem to optimality is exactly the same, in terms of the sets $W^i_\ell$, as to identify $\ell$ such that $\bar{\mu}\in W^i_\ell$ for some $i$. We intend to reduce this task to solving a series of Color Inference problems. We start with presenting the principal building block of our construction -- {\em Individual Inference procedure}.
\paragraph{Individual Inference procedure} is parameterized by $\ell\in\{1,...,L\}$ and a real $\delta>0$.
Given $\ell$ and $\delta$, we initialize the algorithm as follows:
\begin{itemize}
\item mark as red all nonempty sets $W^i_\ell$ along with their elements and the corresponding regular data $(\H,W^i_\ell,\Phi_i\big|_{\H\times W^i_\ell})$;
\item look one by one at all sets $W^i_{\ell'}$ with $i\leq I$ and $\ell'\neq\ell$, and associate with these sets their chunks
\begin{equation}\label{VoronoiChunks}
W^{i\delta}_{\ell\ell'}=\{\mu\in W^i_{\ell'}: u_{\ell\ell'}^T[G\mu]\geq v_{\ell\ell'}+\delta\}.
\end{equation}
Note that the resulting sets are convex and compact. Whenever  $W^{i\delta}_{\ell\ell'}$ is nonempty, we mark  blue this set along with its elements and the corresponding regular data $(\H,W^{i\delta}_{\ell\ell'},\Phi_i\big|_{\H\times W^{i\delta}_{\ell\ell'}})$ .
\end{itemize}
As a result of the above actions, we get a collection of nonempty convex compact subsets $W^{s\delta}_\ell$, $s=1,...,S^\delta_\ell$, of $\bR^d$ and associated regular data ${D^{s\delta}_\ell}=(\H,{W^{s\delta}_\ell},
{\Phi^{s\delta}_\ell})$; the sets ${W^{s\delta}_\ell}$, same as their elements and regular data ${\D^{s\delta}_\ell}$, are colored in red and blue. Note that the sets
 {$W^{s\delta}_\ell$} of different colors
do not intersect (since their images under the mapping $G$ do not intersect), so that a point $\mu\in\bR^m$ gets at most one color. Note also that our collection definitely contains red components.
\par

 Individual Inference Procedure $\A^{\ell\delta}_K$ infers the color of a regular data $\D\in{\mathbf{D}}^{\delta}_\ell=\{\D^{s\delta}_\ell,\,1\leq s\leq S^\delta_\ell\}$  given i.i.d. $K$-repeated observation $\omega^K$ drawn from  a distribution $P\in\S[\D]$:  when the collection ${\mathbf{D}}^{\delta}_\ell$ contains both red and blue regular data, $\A^{\ell\delta}_K$ is exactly Color Inference procedure from section \ref{sect:IC} associated with this collection and our coloring\footnote{The procedure is well defined, since by our assumption, all convex-concave functions we need to deal with are continuous, coercive in the minimization variable, and have closed convex minimization and compact convex maximization domains, so that the required saddle points do exist.}; if no blue regular data is present, $\A^{\ell\delta}_K$
always infers that the color is red.
\par
Observe that if the collection ${\mathbf{D}}^\delta_\ell$ of regular data we have built contains no blue data for some value $\bar{\delta}$ of $\delta$, the same holds true for all $\delta\geq \bar{\delta}$.
Let us define the risk $\widehat{\epsilon}(\ell,\delta)$ of Individual Inference Procedure with parameters $\ell,\delta$ as follows: when $\delta$ is such that ${\mathbf{D}}^\delta_\ell$ contains no blue regular data, $\widehat{\epsilon}(\ell,\delta)=0$, otherwise $\widehat{\epsilon}(\ell,\delta)$ is as stated in Proposition \ref{prop3}. Note that whenever $\delta>0$, $\ell\leq L$, $s\leq S^\delta_\ell$, $\mu\in W^{s\delta}_\ell$ and a probability distribution $P$ satisfies
$$
\bE_{\omega\sim P}\{\exp\{h^T\omega\}\}\leq \Phi^s(h,\mu)\,\,\forall h\in\H,
$$
the quantity
 $\widehat{\epsilon}(\ell,\delta)$ , by construction, is an upper bound on $P$-probability of the event ``as applied to observation $\omega^K=(\omega_1,...,\omega_K)$ with $\omega_1,...,\omega_K$ drawn, independently of each other, from $P$,  $\A^{\ell\delta}_K$ does not recover correctly the color of $\mu$.''
Observe that $\widehat{\epsilon}(\ell,\delta)=0$ for large enough values of $\delta$ (since for large $\delta$ the collection ${\mathbf{D}}^\delta_\ell$ contains no blue data; recall that the parameter sets $\M_i$ are bounded). Besides this, we claim that $\widehat{\epsilon}(\ell,\delta)$ is nonincreasing in $\delta>0$.
\begin{quote}{\small
To support the claim, assume that $\delta'\geq \delta''>0$ are such that {${\mathbf{D}}^{\delta''}_\ell$} contains blue data, and let us show that  $\epsilon':=\widehat{\epsilon}(\ell,\delta')\leq\epsilon'':=\widehat{\epsilon}(\ell,\delta'')$. Indeed, recall that $\epsilon'$ and $\epsilon''$ are $\|\cdot\|_{2,2}$-norms of respective symmetric entry-wise nonnegative matrices $E',E''$ (see Proposition \ref{prop3}). When increasing $\delta$ from $\delta=\delta''$ to $\delta=\delta'$, we reduce the associated $E$-matrix to its submatrix\footnote{The sets $W^{i\delta}_{\ell\ell'}$ shrink as $\delta$ grows, thus  some of the blue sets $W^{i\delta}_{\ell\ell'}$ which are nonempty at $\delta=\delta''$ can become empty when $\delta$ increases to $\delta'$.} and further reduce the entries in this submatrix\footnote{Indeed, by Proposition \ref{prop3}, these entries are obtained from saddle point values of some convex-concave functions by a monotone transformation; it is immediately seen that as $\delta$ grows, these functions and the domains of the minimization argument remain intact, while the domains of the maximization argument shrink, so that the saddle point values cannot increase.}, and thus reduce the norm of the matrix.
}\end{quote}
\paragraph{Aggregation procedure} we propose is as follows:
\begin{quote}
Given tolerance $\epsilon\in (0,\half)$, for every $\ell=1,...,L$ we specify $\delta_\ell>0$, the smaller the better, in such a way that $\widehat{\epsilon}(\ell,\delta_\ell)\leq\epsilon/L$ \footnote{For instance, we could start with  $\delta=\delta^0$ large enough to ensure that $\widehat{\epsilon}(\ell,\delta^0)=0$, and select $\delta_\ell$ as either the last term in the progression $\delta^i={\kappa^i}\delta_0$, $i=0,1,...$, {for some $\kappa\in (0,1)$,} such that $\widehat{\epsilon}(\ell,\delta^i)\leq\epsilon/L$, or the first term in this progression which is ``negligibly small'' (say, less that $10^{-6}$), depending on what happens first.}. Given observation $\omega^K$, we run the procedures $\A^{\ell\delta_\ell}_K$, $1\leq \ell\leq L$. Whenever $ \A^{\ell\delta_\ell}_K$ returns a color, we assign it to the index $\ell$ and to the vector $g_\ell$, so that after all $\ell$'s are processed, some $g_\ell$'s get color ``red,'' some get color ``blue'', and some do not get color. The aggregation procedure returns, as a solution {$\hat{g}(\omega^K)$}, a (whatever)  red vector if one was discovered, and returns, say, $g_1$ otherwise.
\end{quote}
\begin{proposition}\label{prop555}
In the situation and under the assumptions described in the beginning of section \ref{sect:aggr},  let $\omega^K=(\omega_1,...,\omega_K)$ be $K$-element i.i.d. sample drawn from a probability distribution ${\bar{P}}$ which, taken along with some $\bar{i}\leq I$ and $\bar{\mu}\in \M_{\bar{i}}$, satisfies {\rm (\ref{eq1939})}. Then
the ${\bar{P}}$-probability of the event
\begin{equation}\label{aggregationquality}
\|G\bar{\mu}-{\hat{g}(\omega^K)}\|_2\leq \min_{\ell\leq L}\|G\bar{\mu}-g_\ell\|_2+2\max_\ell\delta_\ell
\end{equation}
is at least $1-\epsilon$.
\end{proposition}
 \paragraph{Simple illustration.} Let $I=1$,  $\H=\Omega=\bR^d$, and let $\M=\M_1$ be a nonempty convex compact subset of $\bR^d$. Further, suppose that $\Phi(h;\mu):=\Phi_1(h;\mu)=h^T\mu+{1\over 2}h^T{\Theta}h:\,\H\times \M\to\bR$, where $\Theta$ is a given positive definite matrix. We are also given a $K$-element i.i.d. sample $\omega^K$ drawn from a sub-Gaussian distribution $P$ with sub-Gaussianity parameters $(\mu,\Theta)$. Let also $G\mu\equiv\mu$, so that the aggregation problem we are interested in reads: {\em ``given $L$ estimates $g_1,...,g_L$ of the expectation $\mu$ of a sub-Gaussian random vector $\omega$ with sub-Gaussianity parameters $(\mu,\Theta)$, with a known $\Theta\succ0$,  and $K$-repeated i.i.d. sample $\omega^K$ from the distribution of the vector, we want to select $g_\ell$ which is $\|\cdot\|_2$-closest to the true expectation $\mu$ of $\omega$.''} From now on we assume that $g_1,...,g_L\in \M$ (otherwise, projecting the estimates onto $\M$, we could provably improve their quality) and that $g_1,...,g_L$ are distinct from each other.
 \par
 In our situation, the sets (\ref{Voronoi}), (\ref{VoronoiW}), (\ref{VoronoiChunks}) and functions $\Phi^i$ become
 $$
 \begin{array}{rcl}
 W_\ell&=&\{\mu\in\M: u_{\ell\ell'}^T\mu\leq v_{\ell\ell'},\forall (\ell'\leq L:\ell'\neq\ell)\},\\
 &&u_{\ell\ell'}={g_\ell'-g_\ell\over\|g_\ell'-g_\ell\|_2},\;\;v_{\ell\ell'}={1\over 2}[g_\ell+g_\ell'], \,\ell\neq\ell';\\
 W^\delta_{\ell\ell'}&=&\{\mu\in\M: u_{\ell\ell'}^T\mu\geq v_{\ell\ell'}+\delta\},\,1\leq\ell,\ell'\leq L, \ell\neq\ell',\\
 \Phi(h;\mu)&=&h^T\mu+{1\over 2}h^T\Theta h,
 \end{array}
 $$
 (we are in the case of $I=1$ and thus suppress index $i$ in the notation for $W$'s and $\Phi$). Note that Individual Inference Procedure $\A^{\ell\delta}_K$ deals with exactly one red hypothesis, $\HypSs[\bR^d,W_\ell,\Phi]$, and at most $L-1$ blue hypotheses $\HypSs[\bR^d,W^{\delta}_{\ell\ell'},\Phi]$ associated with nonempty sets $W^{\delta}_{\ell\ell'}$ and $\ell'\neq\ell$.

Applying the construction from section \ref{sect:IC}, we arrive at the aggregation routine as follows (below $\ell,\ell'$ vary in $\{1,...,L\}$):
\begin{itemize}
\item We set
\begin{equation}\label{Goldensluger}
\begin{array}{rcl}
\delta_\ell&=&\max\limits_{\ell'\neq\ell}\sqrt{\ln\left({L\sqrt{L-1}\over\epsilon K}\right)u_{\ell\ell'}^T\Theta u_{\ell\ell'}};\\
w_{\ell\ell'}&=&{1\over 2}\left[g_\ell+g_{\ell'}+\delta_\ell u_{\ell\ell'}\right],\,\ell\neq\ell';\\
\psi_{\ell\ell'}^{(K)}(\omega^K)&=&{\delta_\ell\over 2u_{\ell\ell'}^T\Theta
u_{\ell\ell'}}u_{\ell\ell'}^T\left[Kw_{\ell\ell'}-\sum_{t=1}^K\omega_t\right]+{1\over 2}\ln(L-1),\,\ell\neq\ell'.
\end{array}
\end{equation}
\item Given $\omega^K$, for $1\leq \ell\leq L$ we assign vector $g_\ell$ color ``red,'' if $\psi_{\ell\ell'}^{(K)}(\omega^K)>0$ for all $\ell'\neq\ell$, otherwise we do not assign $g_\ell$ any color;
\item If red vectors were found, we output (any) one of them as solution to the aggregation problem; if no red vectors are found, we output, say, $g_1$ as solution.
\end{itemize}
Proposition \ref{prop555} as applied to the situation in question states that whenever $\omega_1,\omega_2,...\omega_K$ are  drawn, independently from each other, from a sub-Gaussian distribution $P$ with parameters $\mu,\Theta$,  then with $P$-probability at least $1-\epsilon$ the result $\hat{g}(\omega^K)$ of the above aggregation routine satisfies the relation
$$
\|\mu-{\hat{g}(\omega^K)}\|_2\leq \min_{1\leq\ell\leq L}\|\mu-g_\ell\|_2+2\max_{\ell}\delta_\ell,
$$
which essentially recovers the classical $\ell_2$ oracle inequality (cf. \cite[Theorem 4]{Golden2009}).

\section{Beyond the scope of affine detectors}\label{sec:quadratic}
\subsection{Lifted detectors}
The tests developed in sections \ref{sect:pairwise} and \ref{sect:multiple} were based on {\sl affine detectors} -- affine functions $\phi(\omega)$ associated with pairs of composite hypotheses
$H_1:P\in\P_1$, $H_2:P\in\P_2$ on the probability distribution $P$ of observation $\omega\in\Omega=\bR^d$. Such detectors were built to satisfy the relations
\begin{equation}\label{eq1000}
\begin{array}{c}
\int_\Omega\exp\{-\phi_*(\omega)\}P(d\omega)\leq\epsilon_\star\,\,\forall P\in \P_1\ \&\
\int_\Omega\exp\{\phi_*(\omega)\}P(d\omega)\leq\epsilon_\star\,\,\forall P\in \P_2,
\end{array}
\end{equation}
with as small $\epsilon_\star$ as possible (cf. (\ref{eq100})), and affinity of $\phi$ is of absolutely no importance here: all constructions in sections \ref{sect:pairwise}, \ref{sect:multiple} were based upon availability of pairwise detectors
$\phi$, {\sl affine or not}, satisfying, for the respective pairs of composite hypotheses,  relations (\ref{eq100}) with some known $\epsilon_\star$. So far, affinity of detectors was utilized only when {\sl building} detectors satisfying (\ref{eq1000}) via the generic scheme presented in section \ref{sect:situ}.
\par
Now note that given a random  observation $\zeta$  taking values in some $\bR^d$ along with a deterministic function $Z(\zeta):\bR^d\to\bR^D$, we can convert an observation $\zeta$ into an observation
$$
\omega=(\zeta,Z(\zeta))\in\bR^d\times \bR^D.
$$
Here $\omega$ is a deterministic transformation of $\zeta$ which ``remembers'' $\zeta$, so that to make statistical inferences from observations $\zeta$ is exactly the same as to make them from observations $\omega$. However, detectors which are affine in $\omega$ can be nonlinear in $\zeta$: for instance, for $Z(\zeta)=\zeta\zeta^T$, affine in $\omega$ detectors are exactly  detectors {\sl quadratic} in $\zeta$. We see that  within the framework of our approach, passing from $\zeta$ to $\omega$ allows to consider a wider family of detectors and thus to arrive at  a wider family of tests. The potential bottleneck here is the necessity to bring the ``augmented'' observations $(\zeta,Z(\zeta))$ into the scope of our setup.
\paragraph{Example: distributions with bounded support.} Consider the case where the  distribution $P$ of our observation $\zeta\in\bR^d$ belongs to a family $\P$
of Borel probability distributions supported on a given bounded set, for the sake of simplicity -- on the unit Euclidean ball $B$ of $\bR^d$. Given a continuous function $Z(\cdot):\bR^d\to\bR^D$, our goal is to cover the family $\P^+$ of distributions $P^+[P]$ of $\omega=(\zeta,Z(\zeta))$ induced by distributions $P\in \P$ of $\zeta$ by a family $\P[\H,\M,\Phi]$ (or $\S[\H,\M,\Phi]$) associated with some regular data, thus making the machinery we have developed so far applicable to the family of distribution $\P^+$. Assuming w.l.o.g. that
$$
\|Z(z)\|_2\leq 1\;\;\forall (z:\|z\|_2\leq1),
$$
observe that for $P\in \P$ the distribution $P^+[P]$ is sub-Gaussian with sub-Gaussianity parameters $(\theta_P,\Theta=2I_{d+D})$, where
$$
\theta_P=\int_{B} (\zeta,Z(\zeta))P(d\zeta).
$$
It follows that
\begin{quote}
{\sl If we can point out a convex compact set $U$ in $\bR^d\times\bR^D$ such that
\begin{equation}\label{UP}
\theta_P\in U\,\,\forall P\in\P,
\end{equation}
then, specifying regular data as
$$
\begin{array}{rcl}
\H&=&\Omega:=\bR^d\times\bR^D,\\
\M&=&U,\\
\Phi(h;\mu)&=&h^T\mu+h^Th:\H\times \M\to\bR,
\end{array}
$$
we ensure that the family $\S[\H,\M,\Phi]$ contains the family of probability distributions $\P^+=\{P^+[P]:P\in\P\}$ on $\Omega$.}
\end{quote}
How useful is this (by itself pretty crude) observation in the context of our approach depends on how much information on $\P$ can be ``captured'' by a properly selected convex compact set $U$ satisfying (\ref{UP}). We are about to consider in more details ``quadratic lifting'' -- the case where $Z(\zeta)=\zeta\zeta^T$.
\subsubsection{Quadratic lifting, Gaussian case}\label{QLift:Gaussian}
Consider the situation where we are given
\begin{itemize}
\item  a nonempty bounded set $U$ in $\bR^m$;
\item a nonempty convex compact subset $\U$ of the positive semidefinite cone $\bS^d_+$;
\item a matrix $\Theta_*\succ0$ such that $\Theta_*\succeq \Theta$ for all $\Theta\in\U$;
\item an affine mapping $u\mapsto \A(u)=A[u;1]:\bR^m\to\Omega=\bR^d$, where $A$ is a given $d\times(m+1)$ matrix.
\end{itemize}
Now, a pair $(u\in U,\Theta\in\U)$ specifies Gaussian  random vector $\zeta\sim \N(\A(u),\Theta)$ and thus specifies a Borel probability distribution
$P[u,\Theta]$ of $(\zeta,\zeta\zeta^T)$. Let $\Q(U,\U)$ be the family of probability distributions on $\Omega=\bR^d\times\bS^d$ stemming in this fashion
from Gaussian distributions with parameters from $U\times\U$. Our goal is to cover the family $\Q(U,\U)$ by a family of the type $\S[\H,\M,\Phi]$, which, as it was already explained, would allow to use the machinery developed so far in order to decide on pairs of composite Gaussian hypotheses
$$
\begin{array}{l}
H_1:\zeta\sim \N(u,\Theta)\hbox{\ with $(u,\Theta)\in U_1\times\U_1$}\\
H_2:\zeta\sim \N(u,\Theta)\hbox{\ with $(u,\Theta)\in U_2\times\U_2$}\\
\end{array}
$$
via tests based on detectors which are {\sl quadratic} in $\zeta$.
\par
It is convenient to represent a linear form on $\Omega=\bR^d\times\bS^d$ as
$$
h^Tz+\half \Tr(HZ),
$$
where $(h,H)\in\bR^d\times\bS^d$ is the ``vector of coefficients'' of the form, and $(z,Z)\in\bR^d\times\bS^d$ is the argument of the form.
\par
We denote by $b=[0;0;...;0;1]\in \bR^{m+1}$ the last basic orth of $\bR^{m+1}$.
We assume  that for some $\delta\geq0$ it holds
\begin{equation}\label{56delta}
\|\Theta^{1/2}\Theta_*^{-1/2}-I\|\leq \delta\,\,\forall \Theta\in\U,
\end{equation}
where $\|\cdot\|$ is the spectral norm. Observe that for every $\Theta\in\U$ we have $0\preceq \Theta_*^{-1/2}\Theta\Theta_*^{-1/2}\preceq I$, whence $\|\Theta^{1/2}\Theta_*^{-1/2}\|\leq 1$, that is, (\ref{56delta}) is always satisfied with $\delta=2$. Thus, we can assume w.l.o.g. that $\delta\in[0,2]$.
Finally, we  set $B=\left[\begin{array}{c}A\cr b^T\cr\end{array}\right]$.\par
A desired ``covering'' of $\Q(U,\U)$ is given by the following
\begin{proposition}\label{propGausslift} In the notation and under the assumptions of this section, let $\gamma\in(0,1)$ and a convex compact computationally tractable set $\Z\subset\{W\in\bS^{m+1}_+:W_{m+1,m+1}=1\}$ be given, and let ${\Z}$ contain all matrices $Z(u):=[u;1][u;1]^T$ with $u\in U$. Denoting by $\phi_{\Z}(\cdot)$ the support function of $\Z$:
$$
\phi_{\Z}(W)=\max_{Z\in\Z}\Tr(ZW):\bS^{m+1}\to\bR,
$$
let us set
\begin{equation}\label{eq132546}
\begin{array}{rcl}
\H&=&\H_\gamma:=\{(h,H)\in\bR^d\times \bS^d: -\gamma\Theta_*^{-1}\preceq H\preceq \gamma\Theta_*^{-1}\}\\
\M&=&\U,\\
\Phi(h,H;\Theta)&=&-\half \ln\Det(I-\Theta_*^{1/2}H\Theta_*^{1/2})+\half \Tr([\Theta-\Theta_*]H)+{\delta(2-\delta)\over 2(1-\gamma)}\|\Theta_*^{1/2}H\Theta_*^{1/2}\|_F^2\\
&&+\Gamma_{\Z}(h,H):\H\times\U\to\bR,\quad\qquad\hfill\hbox{\rm[$\|\cdot\|_F$ is the Frobenius norm]}\\
\Gamma_{\Z}(h,H)&=&\half \phi_{\Z}\left(bh^TA+A^Thb^T+A^THA+B^T[H,h]^T[\Theta_*^{-1}-H]^{-1}[H,h]B\right)
\\
&=&\half\phi_{\Z}\left(B^T\left[\hbox{\small$\left[\begin{array}{c|c}H&h\cr\hline h^T&\end{array}\right]$}+
\left[H,h\right]^T[\Theta_*^{-1}-H]^{-1}\left[H,h\right]\right]B\right).\\
\end{array}
\end{equation}
Then $\H,\M,\Phi$ form a regular data, and
\begin{equation}\label{target123}
\Q(U,\U)\subset\S[\H,\M,\Pi].
\end{equation}
Besides this, function $\Phi(h,H;\Theta)$ is coercive in the convex argument: whenever $(h_i,H_i)\in\H$ and $\|(h_i,H_i)\|\to\infty$ as $i\to\infty$, we have $\Phi(h_i,H_i;\Theta)\to\infty$, $i\to\infty$, for every $\Theta\in\U$.
\end{proposition}
For proof, see Appendix \ref{proof.propGausslift}. Note that every quadratic constraint $u^TQu+2q^Tu+r\leq 0$ which is valid on $U$ induces a {\sl linear constraint} $\Tr\left(\hbox{\footnotesize$\left[\begin{array}{c|c}Q&q\cr\hline q^T&r\cr\end{array}\right]$}Z\right)\leq 0$ which is valid for all matrices $Z(u)$, $u\in U$, and thus can be incorporated into the description of $\Z$.
\paragraph{Special case.} In the situation of Proposition \ref{propGausslift}, let $u$ vary in a convex compact set $U$. In this case, the simplest way to define $\Z$ such that $Z(u)\in \Z$ for all $u\in U$ is to set
$$
\Z=\left\{W=\hbox{\footnotesize$\left[\begin{array}{c|c}V&u\cr\hline u^T&1\cr\end{array}\right]$}:W\succeq0,u\in U\right\}.
$$
Let us compute the function $\Phi(h,0;\Theta)$. Setting $A=[\bar{A},a]$, where $a$ is the last column of $A$, direct computation yields
{\footnotesize
\bse
\Phi(h,0;\Theta)
&=&{1\over 2}\phi_{\Z}\left(\left[\begin{array}{c|c}&\bar{A}^Th\cr\hline h^T\bar{A}&2a^Th+h^T\Theta_*h\cr\end{array}\right]\right)
={1\over 2}\max\limits_{V,u}\left\{2u^T\bar{A}^Th+2a^Th+h^T\Theta_*h:\left[\begin{array}{c|c}V&u\cr\hline u^T&1\cr\end{array}\right]\succeq0,u\in U\right\}\\
&=&\max\limits_{u\in U}\left\{h^TAu+{1\over 2}h^T\Theta_*h\right\}.
\ese}
Now imagine that we are given two collections $(A_\chi,U_\chi,\U_\chi)$, $\chi=1,2$, of the $(A,U,\U)$-data, with the same number of rows in $A_1$ and $A_2$ and
have associated with $\U_\chi$ $\succeq$-upper bounds $\Theta_{*,\chi}$ on matrices $\Theta\in\U_\chi$. We want to use Proposition \ref{propGausslift} to build {\sl an affine} detector capable to decide on the hypotheses $H_1$ and $H_2$ on the distribution of observation $\omega$, with $H_\chi$ stating that this observation is $\N(A_\chi[u;1],\Theta)$ with some $u\in U_\chi$ and $\Theta\in\U\chi$. To this end, we have to solve the convex-concave saddle point problem
$$
{\cal SV}=\min\limits_h\max\limits_{{\Theta_1\in\U_1,\atop\Theta_2\in \U_2}}{1\over 2}\left[\Phi_1(-h,0;\Theta_1)+\Phi_2(h,0;\Theta_2)\right],
$$
where $\Phi_1$, $\Phi_2$ are the functions associated, as explained in Proposition \ref{propGausslift}, with the first, respectively, with the second collection of the $(A,U,\U)$-data. In view of the above computation, this boils down to the necessity to solve the convex minimization problem
$$
{\cal SV}=\min_h\left\{\max\limits_{{u_1\in U_1,\atop u_2\in U_2}}{1\over 2}\left[{1\over 2}h^T\Theta_{*,1}h+{1\over 2}h^T\Theta_{*,2}h+h^T[A_2[u_2;1]-A_1[u_1;1]]\right]\right\}.
$$
An optimal solution $h_*$ to this problem induces the affine detector
$$
\phi_*(\omega)=h_*^T\omega+a,\,\,a={1\over 2}\left[{1\over 2}h_*^T[\Theta_{*,1}-\Theta_{*,2}]h_*+\max_{u_1\in U_1}[-h_*^TA_1[u_1;1]]-\max_{u_2\in U_2}h_*^TA_2[u_2;1]\right],
$$
and the risk of this detector on the pair of families $\G_1$, $\G_2$ of Gaussian distributions in question is $\exp\{{\cal SV}\}$.
\par
On the other hand, we could build affine detector for the families $\G_1$, $\G_2$ by the machinery from section \ref{subGaussian}, that is, by solving convex-concave saddle point problem
$$
\overline{\cal SV}=\min\limits_h\max\limits_{{u_1\in U_1,\Theta_1\in \U_1,\atop
u_2\in U_2,\Theta_2\in \U_2}}{1\over 2}\left[-h^TA_1[u_1;1]+h^TA_2[u_2,1]+{1\over 2}h^T\Theta_1 h +{1\over 2}h^T\Theta_2 h\right];
$$
the risk of the resulting affine detector on $\G_1$, $\G_2$ is $\exp\{\overline{\cal SV}\}$.
Now assume that
\begin{quote}
(!) $\U_\chi$, $\chi=1,2$, have $\succeq$-maximal elements, and these elements are selected as $\Theta_{*,\chi}$.
\end{quote} In this case the above computation says that ${\cal SV}$ and $\overline{\cal SV}$ are the minimal values of identically equal to each other functions and thus are equal to each other. Thus, in the case of (!) the machinery of Proposition \ref{propGausslift} produces a quadratic detector which can be only better, in terms of risk, than the affine detector yielded by Proposition \ref{propsubGauss}.\footnote{This seems to be tautology -- there are more quadratic detectors than affine ones, so that the best risk achievable with quadratic detectors can be only smaller than the best risk achievable with affine detectors. The point, however, is that Proposition \ref{propGausslift} does not guarantee building the best, in terms of its risk, quadratic detector, it deals with ``computationally tractable approximation'' of this problem. As a result, the quadratic detector constructed in the latter proposition can, in principle, be worse than the affine detector yielded by Proposition \ref{propsubGauss}.}

\paragraph{Numerical illustration.} To get an impression of the performance of quadratic detectors as compared to affine ones in the case of (!), we present here the results of experiment where $U_1=U_1^\rho=\{u\in\bR^{12}:u_i\geq\rho,1\leq i\leq 12\}$, $U_2=U_2^\rho=-U_1^\rho$, $A_1=A_2\in\bR^{8\times 12}$, and $\U_\chi=\{\Theta_{*,\chi}=\sigma_\chi^2I_8\}$ are singletons. The risks of affine, quadratic and ``purely quadratic'' (with $h$ set to 0) detectors on the pair $\G_1,\G_2$ of families of Gaussian distributions, with $\G_\chi=\{\N(\theta,\Theta_{*\chi}): \theta\in A_\chi U_\chi^\rho\}$, ar given in Table \ref{t:1}.
\begin{table}
\begin{center}
{
\begin{tabular}{||c||c|c||c|c|c||}
\hline
$\rho$&$\sigma_1$&$\sigma_2$&\shortstack{unrestricted\\ $H$ and $h$}&$H=0$&$h=0$\\
\hline
0.5&2&2&0.31&0.31&1.00\\
\hline
0.5&1&4&0.24&0.39&0.62\\
\hline
0.01&1&4&0.41&1.00 &0.41\\
\hline
\end{tabular}
}
\caption{\label{t:1} Risk of quadratic detector}
\end{center}
\end{table}

We see that
\begin{itemize}
\item when deciding on families of Gaussian distributions with common covariance matrix and expectations varying in associated with the families convex sets, passing from affine detectors described by Proposition \ref{propsubGauss} to quadratic detectors, does not affect the risk (first row in the table). This is a general fact: by the results of \cite{GJN}, in the situation in question affine detectors are optimal in terms of risk among {\sl all possible} detectors.
\item when deciding on families of Gaussian distributions in the case where distributions from different families can have close expectations (third row in the table), affine detectors are useless, while the quadratic ones are not, provided that $\Theta_{*,1}$ differs from $\Theta_{*,2}$. This is how it should be -- we are in the case where the first moments of the distribution of the observation bear no definitive information on the family this distribution belongs to, which makes affine detectors useless. In contrast, quadratic detectors are able to utilize information (valuable when $\Theta_{*,1}\neq\Theta_{*,2})$ ``stored'' in the second moments of the observation.
\item ``in general'' (second row in the table), both affine and purely quadratic components in a quadratic detector are useful; suppressing one of them can increase significantly the attainable risk.
\end{itemize}
\subsubsection{Quadratic lifting: Bounded observations}

It is convenient to represent a ``quadratically lifted observation'' $(\zeta,\zeta\zeta^T)$ by the matrix
$$
Z(\zeta)=\left[\begin{array}{c|c}\zeta\zeta^T&\zeta\cr\hline\zeta^T&1\cr\end{array}\right]\in\bS^{d+1}.
$$
Assume that all distributions from $\P$ are supported on the solution set $\X$ of a system of quadratic constraints
$$
f_\ell(\zeta):=\zeta^TA_\ell\zeta+2a_\ell^T\zeta+\alpha_\ell\leq 0,\,1\leq\ell\leq L,
$$
where $A_\ell\in\bS^{d}$ are such that $\sum_\ell\bar{\lambda}_\ell A_\ell\succ0$ for properly selected $\bar{\lambda}_\ell\geq0$; as a consequence, $\X$ is bounded (since $\bar{f}(\zeta):=\sum_{\ell=1}^L\bar{\lambda}_\ell f_\ell(\zeta)$ is a strongly convex quadratic form which is $\leq0$ on $X$). Setting
$$
Q_0=\left[\begin{array}{c|c}&\cr\hline &1\cr\end{array}\right],Q_\ell=\left[\begin{array}{c|c}A_\ell&a_\ell\cr\hline a_\ell^T&\alpha_\ell\cr\end{array}\right],\,1\leq \ell\leq L,
$$
observe that the distribution $P^+$ of $Z(\zeta)$ induced by a distribution $P\in\P$ is supported on the closed convex set
$$
\X^+=\{Z\in\bS^{d+1}:\; Z\succeq0,\; \Tr(Q_0Z)=1,\;\Tr(Q_\ell Z)\leq0,\;\;\ell=1,...,L\}
$$
which is bounded\footnote{indeed, for large $\lambda_0>0$, the matrix $Q=\lambda_0Q_0+\sum_{\ell=1}^L\bar{\lambda}_\ell Q_\ell$ is positive definite, and we conclude that $\X^+$ is contained in the {\sl bounded} set $\{Z\in\bS^{d+1}:Z\succeq0,\Tr(QZ)\leq\lambda_0\}$.}. The support function of this set is
$$
\phi_{X^+}(H)=\max_{Z\in\X^+}\Tr(HZ)=\max_Z\left\{\Tr(HZ):\left\{\begin{array}{l}
Z\succeq 0,Z_{d+1,d+1}=1\\
\Tr(Q_\ell Z)\leq 0,\,1\leq\ell\leq L\\
\end{array}\right.\right\}
$$
Recalling  Example 4 in section \ref{sect:basicexamples} and section \ref{sect:Support}, we arrive at the regular data
$$
\H=\bS^{d+1},\,\M=\X^+,\,\Phi(h,\mu)=\inf_{g\in \bS^{d+1}}\left\{\Tr((h-g)\mu)+{\eight}[\phi_{X_+}(h-g)+\phi_{X_+}(g-h)]^2+\phi_{\X_+}(g)\right\}
$$
such that
$$\forall P\in\P: P^+\in\S[\bS^{d+1},\{e[P^+]\},\Phi(\cdot,e[P^+])],$$
where $e[P^+]$ is the expectation of $P^+$, and therefore
$$
\P^+=\{P^+:P\in \P\}\subset \S[\bS^{d+1},\M,\Phi].
$$

\appendix
\section{Proofs}
\subsection{Proof of Proposition \ref{propbounded}}\label{app1}
All we need is to verify (\ref{eq090807}) and to check that the right hand side function in this relation is
convex. The latter is evident, since $\phi_X(h)+\phi_X(-h)\geq 2\phi_X(0)=0$ and $\phi_X(h)+\phi_X(-h)$ is convex. To verify (\ref{eq090807}), let us fix $P\in\P[X]$ and $h\in\bR^d$ and set $$\nu=h^Te[P],$$ so that $\nu$ is the expectation of
$h^T\omega$ with $\omega\sim P$. Note that $-\phi_X(-h)\leq\nu\leq\phi_X(h)$, so that (\ref{eq090807}) definitely holds true when $\phi_X(h)+\phi_X(-h)=0$. Now let
$$\eta:=\half \left[\phi_X(h)+\phi_X(-h)\right]>0,$$
and let
$$
a=\half \left[\phi_X(h)-\phi_X(-h)\right],\;\;\beta=(\nu-a)/\eta.
$$
Denoting by $P_h$ the distribution of $h^T\omega$ induced by the distribution $P$ of $\omega$ and noting that this distribution is supported on $[-\phi_X(-h),\phi_X(h)]=[a-\eta,a+\eta]$ and has expectation $\nu$, we get
$$
\beta\in[-1,1]$$
and
$$
\gamma:=\int\exp\{h^T\omega\}P(d\omega)=\int_{a-\eta}^{a+\eta}[\e^s-\lambda(s-\nu)]P_h(ds)
$$
for all $\lambda\in\bR$.
Hence,
\bse
\ln(\gamma)&\leq& \inf_\lambda\ln\left(\max_{a-\eta\leq s\leq a+\eta} [\e^s-\lambda(s-\nu)]\right)=a+\inf_\rho\ln\left(\max_{-\eta\leq t\leq\eta}[\e^t-\rho(t-[\nu-a])]\right)\\
&=&a+\inf_\rho\ln\left(\max_{-\eta\leq t\leq\eta}[\e^t-\rho(t-\eta\beta)]\right)\leq
a+\ln\left(\max_{-\eta\leq t\leq\eta}[\e^t-\bar{\rho}(t-\eta\beta)\right)
\ese
with $\bar{\rho}=(2\eta)^{-1}(\e^\eta-\e^{-\eta})$. 
The function $g(t)=\e^t-\bar{\rho}(t-\eta\beta)$ is convex on $[-\eta,\eta]$, and
\[g(-\eta)=g(\eta)=
\cosh(\eta)+\beta\sinh(\eta),
\] which combines with the above computation to yield the relation
\begin{equation}\label{eq987098}
\ln(\gamma)\leq a+\ln(\cosh(\eta)+\beta\sinh(\eta)),
\end{equation}
and all we need to verify is that
\begin{equation}\label{eq908070}
 \forall (\eta>0,\beta\in[-1,1]): \;\beta\eta +\half \eta^2-\ln(\cosh(\eta)+\beta\sinh(\eta))\geq 0.
\end{equation}
Indeed, if \rf{eq908070} holds true (\ref{eq987098}) implies that
$$
\ln(\gamma) \leq a+\beta\eta+\half \eta^2=\nu+\half \eta^2,$$
 which, recalling what $\gamma$, $\nu$ and $\eta$ are, is exactly what we want to prove.
 \par Verification of (\ref{eq908070}) is as follows.
The left hand side in (\ref{eq908070}) is convex in $\beta$ for $\beta>-{\cosh(\eta)\over \sinh(\eta)}$ containing, due to $\eta>0$, the range of $\beta$ in (\ref{eq908070}). Furthermore, the minimum of the left hand side of (\ref{eq908070}) over $\beta>-\coth(\eta)$ 
is attained when $\beta={\sinh(\eta)-\eta\cosh(\eta)\over\eta\sinh(\eta)}$ and is equal to
$$
r(\eta)=\half\eta^2+1-\eta\coth(\eta)
-\ln(\sinh(\eta)/\eta).
$$ All we need to prove is that the latter quantity is nonnegative whenever $\eta>0$. We have
$$
r'(\eta)=\eta-\coth(\eta)
-\eta(1-\coth^2(\eta))
-\coth(\eta)
+\eta^{-1}
=(\eta\coth(\eta)
-1)^2\eta^{-1}\geq0,
$$
and since $r(+0)=0$, we get $r(\eta)\geq0$ when $\eta>0$. \hfill$\Box$
\subsection{Proof of Proposition \protect{\ref{prop555}}}\label{app:prop555}
Let ${\mathbf{D}}_\ell$ denote the collection of regular data processed by $\A^{\ell\delta_\ell}_K.$ Let, further, the $K$-repeated random observation $\omega^K$, the probability distribution ${\bar{P}}$, index $\bar{i}\leq I$ and vector $\bar{\mu}\in\M_{\bar{i}}$ be as in the premise of the proposition,  let $\bar{g}=G\bar{\mu}$, and let $\ell_*\leq L$ be such that
$\|\bar{g}-g_{\ell_*}\|_2\leq \|\bar{g}-g_\ell\|_2$, $\ell=1,...,L$. Finally, let $\E$ be the event ``for every $\ell=1,...,L$ such that {$\bar{P}$} obeys one or more of the hypotheses $\HypSs[\D]$, $\D\in{\mathbf{D}}_\ell$, processed by procedure $\A^{\ell\delta_{\ell}}_K$, this procedure correctly recovers the color of these hypotheses.'' By construction and due to the union bound, the ${\bar{P}^K}$-probability of $\E$ is at least $1-\epsilon$. It follows that all we need to verify the claim of the proposition is to show that when $\omega^K\in\E$, relation (\ref{aggregationquality}) takes place. Thus, let us fix $\omega^K\in\E$.
\par
Observe, first, that $\bar{g}\in V_{\ell_*}$, whence $\bar{\mu}\in W^{\bar{i}}_{\ell_*}$. Thus, when running
$\A^{\ell_*\delta_{\ell_*}}_K$, ${\bar{P}^K}$ obeys a red one among the hypotheses $\HypSs[\D]$, $\D\in{\mathbf{D}}_{\ell_*}$ processed by $\A^{\ell_*\delta_{\ell_*}}_K$, and since we are in the case of $\omega^K\in \E$, $g_{\ell_*}$ gets a color, namely, color ``red.'' By construction of our aggregation procedure, its output can be either $g_{\ell_*}$ -- and in this case (\ref{aggregationquality}) clearly holds true, or another vector, let it be denoted $g_{\ell_+}$ ($\ell_+\neq\ell_*$), which was also assigned red color. We claim that the vector $\bar{g}=G\bar{\mu}$ satisfies the relation
\begin{equation}\label{condition2}
u_{\ell_+\ell_*}^T\bar{g}<v_{\ell_+\ell_*}+\delta_{\ell_+}.
\end{equation}
Indeed, otherwise we have $\bar{\mu}\in W^{\bar{i}\delta_{\ell_+}}_{\ell_+\ell_*}$, meaning that ${\bar{P}}^K$ obeys a hypothesis $\HypSs[\D]$ processed when running $\A^{\ell_+\delta_{\ell_+}}_K$ (i.e., with $\D\in{\mathbf{D}}_{\ell_+}$), and this hypothesis is blue. Since we are in the case of $\E$, this implies that the color inferred by $\A^{\ell_+\delta_{\ell_+}}_K$ is ``blue,'' which is a desired contradiction.
\par
Now we are nearly done: indeed, $g_{\ell_*}$ is the $\|\cdot\|_2$ closest to $\bar{g}$ point among $g_1,...,g_L$, implying that
\begin{equation}\label{condition3}
u_{\ell_+\ell_*}^T\bar{g}\geq v_{\ell_+\ell_*}.
\end{equation}
Recalling what $u_{\ell_*\ell_+}$ and $v_{\ell_*\ell_+}$ are, the relations (\ref{condition2}) and (\ref{condition3}) tell us the following story about the points $\bar{g}$,
 $g_*:=g_{\ell_*}$, $g_+:=g_{\ell_+}$ and the hyperplane $H=\{g\in \bR^m:\,u_{\ell_*\ell_+}g=v_{\ell_*\ell_+}$\}: $g_*$ and $g_+$ are symmetric to each other w.r.t. $H$, and $\bar{g}$ is at most at the distance $\delta:=\delta_{\ell_+}$ of $H$. An immediate observation is that in this case
 \be
 \|g_+-\bar{g}\|_2\leq \|g_*-\bar{g}\|_2+2\delta,
 \ee{eq:arik100}
 and we arrive at (\ref{aggregationquality}). \par
 To justify {\rf{eq:arik100}} note that by shift and rotation we can reduce the situation to the one when $g_*,g_+,\bar{g}$ belong to the linear span of the first two basic orhts {with} the first two coordinates of these three vectors {being}, respectively, $[-r;0]$ (for $g_*$), $[r;0]$ (for $g_+$) and
 $
 [d;h]$ (for $\bar{g}$), with $|d|\leq\delta$. Hence
 $$
 \begin{array}{l}
 \|g_+-\bar{g}\|_2-\|g_*-\bar{g}\|_2={{[(r-d)^2+h^2]-[(r+d)^2+h^2]}\over \|g_+-\bar{g}\|_2+\|g_*-\bar{g}\|_2}=
 {-4rd\over \|g_+-\bar{g}\|_2+\|g_*-\bar{g}\|_2}\leq {4r\delta\over\|g_+-g_*\|_2}
 =2\delta,
 \end{array}
 $$
 as claimed. {\hfill $\Box$}

\subsection{Proof of Proposition \ref{propGausslift}}\label{proof.propGausslift}
\paragraph{1$^0$.} For any $u\in\bR^m$, $h\in\bR^d,\;\Theta\in\bS^d_+$ and $H\in\bS^d$ such that  $-I\prec \Theta^{1/2}H\Theta^{1/2}\prec I$ we have
\be
\lefteqn{\Psi(h,H;u,\Theta):=\ln\left(\bE_{\zeta\sim\N(\A(u),\Theta)}\left\{\exp\{h^T\zeta+\half \zeta^TH\zeta\}\right\}\right)}\nn
&=&\ln\left(\bE_{\xi\sim\N(0,I)}\left\{\exp\{h^T[\A(u)+\Theta^{1/2}\xi]+\half [\A(u)+\Theta^{1/2}\xi]^TH[\A(u)+\Theta^{1/2}\xi]\right\}\right) \nn
&=&-\half \ln\Det(I-\Theta^{1/2}H\Theta^{1/2})\nn
&&+h^T\A(u)+\half \A(u)^TH\A(u)+\half [H\A(u)+h]^T\Theta^{1/2}[I-\Theta^{1/2}H\Theta^{1/2}]^{-1}\Theta^{1/2}[H\A(u)+h]
\nn
&=&-\half \ln\Det(I-\Theta^{1/2}H\Theta^{1/2})
+\half [u;1]^T\left[bh^TA+A^Thb^T+A^THA\right][u;1]\nn
&&+\half [u;1]^T\left[B^T[H,h]^T\Theta^{1/2}[I-\Theta^{1/2}H\Theta^{1/2}]^{-1}\Theta^{1/2}[H,h]B\right][u;1]
\ee{56eq2}
(because $h^T\A(u)=[u;1]^Tbh^TA[u;1]=[u;1]^TA^Thb^T[u;1]$ and $H\A(u)+h=[H,h]B[u;1]$).
\par
 Observe that when $(h,H)\in \H_\gamma$, we have $\Theta^{1/2}[I-\Theta^{1/2}H\Theta^{1/2}]^{-1}\Theta^{1/2}=[\Theta^{-1}-H]^{-1}\preceq
[\Theta_*^{-1}-H]^{-1},$ so that (\ref{56eq2}) implies that for all $u\in\bR^m,\;\Theta\in\U,$ and $(h,H)\in\H_\gamma$,
\be
\Psi(h,H;u,\Theta)&\leq&-\half \ln\Det(I-\Theta^{1/2}H\Theta^{1/2})\nn
&&+\half [u;1]^T\underbrace{\left[bh^TA+A^Thb^T+A^THA+B^T[H,h]^T[\Theta_*^{-1}-H]^{-1}[H,h]B\right]}_{Q[H,h]}[u;1]\nn
&=&-\half \ln\Det(I-\Theta^{1/2}H\Theta^{1/2})+\half \Tr(Q[H,h]Z(u))\nn
&\leq&-\half \ln\Det(I-\Theta^{1/2}H\Theta^{1/2})+\Gamma_{\Z}(H,h)
\ee{56eq22}
(we have taken into account that $Z(u)\in\Z$ when $u\in U$ (premise of the proposition) and therefore $\Tr(Q[H,h]Z(u))\leq\phi_{\Z}(Q[H,h])$).
\paragraph{2$^0$.} Now let us upper-bound the function
$$G(h,H;\Theta)=-\half \ln\Det(I-\Theta^{1/2}H\Theta^{1/2})$$
on the domain $(h,H)\in\H_\gamma$, $\Theta\in\U$. For $(h,H)\in\H_\gamma$ and $\Theta\in\U$ fixed we have
\be
\|\Theta^{1/2}H\Theta^{1/2}\|&=&\|[\Theta^{1/2}\Theta_*^{-1/2}][\Theta_*^{1/2}H\Theta_*^{1/2}][\Theta^{1/2}\Theta_*^{-1/2}]^T\|\nn
&\leq& \|\Theta^{1/2}\Theta_*^{-1/2}\|^2\|\Theta_*^{1/2}H\Theta_*^{1/2}\|\leq \|\Theta_*^{1/2}H\Theta_*^{1/2}\|=:d(H)\leq \gamma
\ee{56eq1234}
(we have used the fact that $0\preceq [\Theta^{1/2}\Theta_*^{-1/2}]^T[\Theta^{1/2}\Theta_*^{-1/2}]\preceq I$ due to $0\preceq \Theta\preceq \Theta_*$, whence $\|\Theta^{1/2}\Theta_*^{-1/2}\|\leq1$).
Denoting by $\|\cdot\|_F$ the Frobenius norm of a matrix and noting that $\|AB\|_F\leq\|A\|\|B\|_F$, computation completely similar to the one in (\ref{56eq1234}) yields
\begin{equation}\label{56eq1235}
\|\Theta^{1/2}H\Theta^{1/2}\|_F\leq \|\Theta_*^{1/2}H\Theta_*^{1/2}\|_F=:D(H).
\end{equation}
Besides this, setting $F(X)=-\ln\Det(X):\inter\bS^d_+\to\bR$ and equipping $\bS^d$ with the Frobenius inner product, we have $\nabla F(X)=-X^{-1}$, so that with $R_0=\Theta_*^{1/2}H\Theta_*^{1/2}$, $R_1=\Theta^{1/2}H\Theta^{1/2}$, and $\Delta=R_1-R_0$, we have for properly selected $\lambda\in(0,1)$ and $R_\lambda=\lambda R_0+(1-\lambda)R_1$:
\bse
F(I-R_1)&=&F(I-R_0-\Delta)=F(I-R_0)+\langle \nabla F(I-R_\lambda),-\Delta\rangle =F(I-R_0)+\langle(I-R_\lambda)^{-1},\Delta\rangle\\
&=&F(I-R_0)+\langle I,\Delta\rangle +\langle (I-R_\lambda)^{-1}-I,\Delta\rangle.
\ese
We conclude that
\begin{equation}\label{56eqtru}
F(I-R_1)\leq F(I-R_0)+\Tr(\Delta)+\|I-(I-R_\lambda)^{-1}\|_F\|\Delta\|_F.
\end{equation}
Denoting by $\mu_i$ the eigenvalues of $R_\lambda$ and noting that $\|R_\lambda\|\leq \max[\|R_0\|,\|R_1\|]=d(H)\leq\gamma$ (see (\ref{56eq1234})), we have $|\mu_i|\leq\gamma$, and therefore eigenvalues $\nu_i=1-{1\over 1-\mu_i}=-{\mu_i\over 1-\mu_i}$ of
$I-(I-R_\lambda)^{-1}$ satisfy $|\nu_i|\leq |\mu_i|/(1-\mu_i)\leq |\mu_i|/(1-\gamma)$, whence $$\|I-(I-R_\lambda)^{-1}\|_F\leq \|R_\lambda\|_F/(1-\gamma).$$
Noting that $\|R_\lambda\|_F\leq\max[\|R_0\|_F,\|R_1\|_F]\leq D(H)$, see (\ref{56eq1235}), we conclude that $\|I-(I-R_\lambda)^{-1}\|_F\leq D(H)/(1-\gamma)$, so that (\ref{56eqtru}) yields
\begin{equation}\label{56eqtru1}
F(I-R_1)\leq F(I-R_0)+\Tr(\Delta)+D(H)\|\Delta\|_F/(1-\gamma).
\end{equation}
Further, by (\ref{56delta}) the matrix $D =\Theta^{1/2}\Theta_*^{-1/2}-I$ satisfies $\|D\|\leq\delta$, whence
$$
\Delta=\underbrace{\Theta^{1/2}H\Theta^{1/2}}_{R_1}-\underbrace{\Theta_*^{1/2}H\Theta_*^{1/2}}_{R_0}=
(I+D)R_0(I+D^T)-R_0=DR_0+R_0D^T+DR_0D^T.
$$
Consequently,
$$
\|\Delta\|_F\leq \|DR_0\|_F+\|R_0D^T\|_F+\|DR_0D^T\|_F\leq [2\|D\|+\|D\|^2]\|R_0\|_F\leq \delta(2+\delta)\|R_0\|_F=
\delta(2+\delta)D(H).
$$
This combines with (\ref{56eqtru1}) and the relation
\[
\Tr(\Delta)=\Tr(\Theta^{1/2}H\Theta^{1/2}-\Theta_*^{1/2}H\Theta_*^{1/2})=\Tr([\Theta-\Theta_*]H)
\] to yield
$$
F(I-R_1)\leq F(I-R_0)+\Tr([\Theta-\Theta_*]H)+{\delta(2+\delta)\over1-\gamma}\|\Theta_*^{1/2}H\Theta_*^{1/2}\|_F^2,
$$
and we conclude that for all $(h,H)\in\H_\gamma$ and $\Theta\in\U$,
\be
\lefteqn{G(h,H;\Theta)=\half F(I-R_1)}\nn
&\leq&\widehat{G}(h,H;\Theta):= -\half \ln\Det(I-\Theta_*^{1/2}H\Theta_*^{1/2})
+\half \Tr([\Theta-\Theta_*]H)+{\delta(2+\delta)\over 2(1-\gamma)}\|\Theta_*^{1/2}H\Theta_*^{1/2}\|_F^2.
\ee{56arrive11}
Note that $\widehat{G}(h,H;\Theta)$ clearly is convex-concave and continuous on $\H\times\M=\H_\gamma\times\U$.
\paragraph{3$^0$.} Combining (\ref{56arrive11}), (\ref{56eq22}), (\ref{eq132546}) and the origin of $\Psi$, see (\ref{56eq2}), we arrive at
$$
\forall ((u,\Theta)\in U\times\U,(h,H)\in\H_\gamma=\H):\;\;
\ln\left(\bE_{\zeta\sim\N(u,\Theta)}\left\{\exp\{h^T\zeta+\half \zeta^TH\zeta\}\right\}\right)\leq\Phi(h,H;\Theta).
$$
This is all we need, up to verification of the claim that $\H,\M,\Phi$ is regular data, which boils down to checking that $\Phi:\H\times\M\to\bR$ is convex-concave and continuous. The latter check, recalling that $\widehat{G}(h,H;\Theta):\H\times\M$ indeed is convex-concave and continuous, reduces to verifying that $\Gamma(h,H)$ is convex and continuous on $\H_\gamma$. Recalling that $\Z$ is nonempty compact set, the function $\phi_{\Z}(\cdot):\bS^{d+1}$ is continuous, implying the continuity of $\Gamma(h,H)=\half\phi_{\Z}(Q[H,h])$ on $\H=\H_\gamma$ ($Q[H,h]$ is defined in (\ref{56eq22})). To prove convexity of $\Gamma$, note that $\Z$ is contained in $\bS^{m+1}_+$, implying that $\phi_{\Z}(\cdot)$ is convex and $\succeq$-monotone. On the other hand, by Schur Complement Lemma, we have
$$
\begin{array}{rcl}
S&:=&\{(h,H,G): G\succeq Q[H,h],(h,H)\in\H_\gamma\}\\
&=&\left\{(h,H,G):\left[\begin{array}{c|c}G-[bh^TA+A^ThB+A^THA]&B^T[H,h]^T\cr\hline
[H,h]B&\Theta_*^{-1}-H\cr\end{array}\right]\succeq0,(h,H)\in\H_\gamma\right\},\\
\end{array}
$$
implying that $S$ is convex. Since $\Phi_{\Z}(\cdot)$ is $\succeq$-monotone, we have
$$
\{(h,H,\tau):(h,H)\in \H_\gamma,\;\tau\geq\Gamma(h,H)\}=\{(h,H,\tau):\;
\exists G: G\succeq Q[H,h],\;2\tau\geq \phi_{\Z}(G),\;(h,H)\in\H_\gamma\},
$$
and we see that the epigraph of $\Gamma$ is convex (since the set $S$ and the epigraph of $\phi_{\Z}$ are so), as claimed.
\paragraph{4$^0$.} It remains to prove that $\Phi$ is coercive in $H,h$. Let $\Theta\in\U$ and $(h_i,H_i)\in\H_\gamma$ with $\|(h_i,H_i)\|\to\infty$ as $i\to\infty$, and let us prove that $\Phi(h_i,H_i;\Theta)\to\infty$. Looking at the expression for $\Phi(h_i,H_i;\Theta)$, it is immediately seen that all terms in this expression, except for the terms coming from $\phi_{\Z}(\cdot)$, remain bounded as $i$ grows, so that all we need to verify  is that the $\phi_{\Z}(\cdot)$-term goes to $\infty$ as $i\to\infty$. Observe that $H_i$ are uniformly bounded due to $(h_i,H_i)\in \H_\gamma$, implying that $\|h_i\|_2\to\infty$ as $i\to\infty$. Denoting by $e$ the last basic orth of $\bR^{d+1}$ and by $b$, as before, the last basic orth of $\bR^{m+1}$, note that, by construction, $B^Te=b$. Now let $W\in\Z$, so that $W_{m+1,m+1}=1$. Taking into account that the matrices $[\Theta_*^{-1}-H_i]^{-1}$ satisfy $\alpha I_d\preceq [\Theta_*^{-1}-H_i]^{-1}\preceq \beta I_d$  for some positive $\alpha,\beta$ due to $H_i\in\H_\gamma$, observe that
$$
\begin{array}{l}
\underbrace{\left[\left[\begin{array}{c|c}H_i&h_i\cr\hline h_i^T&\end{array}\right]+
\left[H_i,h_i\right]^T[\Theta_*^{-1}-H_i]^{-1}\left[H_i,h_i\right]\right]}_{Q_i}
=
\underbrace{\left[h_i^T[\Theta_*^{-1}-H_i]^{-1}h_i\right]}_{\alpha_i\|h_i\|_2^2}ee^T+R_i,\\
\end{array}
$$
where $\alpha_i\geq\alpha>0$ and $\|R_i\|_F\leq C(1+\|h_i\|_2)$. As a result,
\bse
\phi_{\Z}(B^TQ_iB)&\geq&\Tr(WB^TQ_iB)=\Tr(WB^T[\alpha_i\|h_i\|_2^2ee^T+R_i]B)\\
&=&\alpha_i\|h_i\|_2^2
\underbrace{\Tr(Wbb^T)}_{=W_{m+1,m+1}=1}-\|BWB^T\|_F\|R_i\|_F
\geq
\alpha\|h_i\|_2^2-C(1+\|h_i\|_2)\|BWB^T\|_F,\\
\ese
and the concluding quantity tends to $\infty$ as $i\to\infty$ due to $\|h_i\|_2\to\infty$, $i\to\infty$.
\hfill $\Box$

\end{document}